\documentclass[submission,copyright,creativecommons]{eptcs}
\providecommand{\event}{ThEdu'22} 

\usepackage{iftex}
\usepackage{epsfig,subfigure}
\usepackage{amsmath,amssymb}
\usepackage{enumerate}

\newtheorem{theorem}{Theorem}

\newtheorem{definition}[theorem]{Definition}

\newtheorem{remark}[theorem]{Remark}

\title{Computer Assisted Proofs and Automated Methods\\ in Mathematics Education}
\author{Thierry Noah Dana-Picard
\institute{Jerusalem College of Technology\\ Jerusalem, Israel}
\email{ndp@jct.ac.il}
}
\def\titlerunning{Computer assisted proofs in Mathematics Education}
\def\authorrunning{Th. Dana-Picard}
\begin{document}
\maketitle

\begin{abstract}
  This survey paper is an expanded version of an invited keynote at
  the ThEdu'22 workshop, August 2022, in Haifa (Israel). After a short
  introduction on the developments of CAS, DGS and other useful
  technologies, we show implications in Mathematics Education, and in
  the broader frame of STEAM Education.  In particular, we discuss the
  transformation of Mathematics Education into
  exploration-discovery-conjecture-proof scheme, avoiding usage as a
  black box . This scheme fits well into the so-called 4 C's of 21st
  Century Education. Communication and Collaboration are emphasized
  not only between humans, but also between machines, and between man
  and machine. Specific characteristics of the outputs enhance the
  need of Critical Thinking. The usage of automated commands for
  exploration and discovery is discussed, with mention of limitations
  where they exist. We illustrate the topic with examples from
  parametric integrals (describing a ``cognitive neighborhood'' of a
  mathematical notion), plane geometry, and the study of plane curves
  (envelopes, isoptic curves). Some of the examples are fully worked
  out, others are explained and references are given.
\end{abstract}

\section{Introduction}
\label{intro}

\subsection{A long history made very short, with a personal touch.}
During the 1980's, a computer-assisted proof was not always easily accepted by the community. At the beginning of the era of machine-assisted computations, only numerical algorithms existed, allowing approximate results. These algorithms could not always provide all the solutions to a given problem.

Later, algorithms for symbolic computations began to be developed and
implemented. Sometimes the researcher had to write in a general
language a specific program for his/her problem. M. Schaps did in the
80's for classification of generic associative
algebras~\cite{dps1993,dps1996} and wrote programs for computing
Hochschild cohomology groups and groups of automorphisms of given
algebras, each one allowing checking the validity of the results of
the other one.  Vanishing of the homology group ensures that the given
algebra is rigid, whence defines a component in the variety
$\textbf{Alg}_n$ which parameterizes the $n-$dimensional associative
algebras. For local algebras, i.e. algebras of the form
$\mathbb{F}[X]/\mathcal{I}$, where $\mathcal{I}$ is an ideal in the
polynomial ring $\mathbb{F}[X]$, we began using computations of
Gr\"obner bases. The project yielded the classification of generic
algebras in dimensions 6, 7 and 8; see~\cite{dps1993}. Smaller
dimensions were studied previously.

Other researchers in the world were already working using computer computations, but at that time a disclaimer was sometimes added about the non-responsibility of the journal regarding the computations. The computations were time-consuming (for one local algebra, the computation could take 20 minutes, with a 4.77Hz PCs. In order to read the output we had to print it (on wide continuous sheets) and we stored the printouts in case somebody would request to check them.

Development and usage of software did not remain the exclusive
property of researchers. Of course, as long as computers were too big
and too expensive, they could not enter the regular classroom. With
the first hand-held devices and the personal computers, things
changed. Hand-held calculators worked numerically, but quite quickly
algorithms for symbolic computations were developed, and implemented
in a variety of devices. The Derive software was contained on a small
floppy disk, a consequence of the choice of algorithms. For example,
the computation by Derive of definite integrals is based on a theorem,
which is both easy to prove and rarely presented in
textbooks~\cite{dp-integrals}. In~\cite{humanvscomputer}, we proposed
to Derive a parametric improper integral. Derive computed it
immediately with the general parameter; at that time, no other system
could do it, but for small integer values of the parameter. Later a
``clone'' of the software was implemented in the TI92. Today, GeoGebra
has versions for PCs, iPads and smartphones. In the recent period,
mobile versions of the different mathematical software allowed the
development of outdoor mathematical
activities~\cite{MathCityTrails,MathTrails}.

Classical methods of teaching mathematics led sometimes teachers to
think that mathematics are a fixed domain of knowledge. Even teacher
trainers eventually claim this. But the world has changed profoundly,
and the various available technologies transformed mathematics into an
experimental domain, and discovering novelties at an early stage of
education.  P. Quaresma says (\cite{quaresma2020}): ``Scientific
research and education at all levels are concerned primarily with the
discovery, verification, communication, and application of scientific
knowledge. Learning, reusing, inventing, and archiving are the four
essential aspects of knowledge accumulation in mankind's civilization
process. In this cycle of knowledge accumulation, which has been
supported for thousands of years by written books and other physical
means, rigorous reasoning has always played an essential role.
Nowadays this process is becoming more and more effective due to the
availability of new paradigms based on computer
applications. Geometric reasoning with such computer applications is
one of the most attractive challenges for future accumulation and
dissemination of knowledge.''

In this paper, we relate mostly to the two kinds of mathematical software, namely Computer Algebra Systems (CAS) and Dynamic Geometry Systems (DGS)\footnote{There exist other kinds of software, proof assistants and theorem provers, whose importance  cannot be overemphasized. We refer the interested reader to recent contributions to ThEdu, also in this volume}.   Actually, the distinction between them fades more and more. For example, GeoGebra began as a DGS, but developed in other directions, and a CAS called Giac~\cite{kovacs-parisse} is embedded in it. New features include also tools for Augmented Reality, allowing outdoor activities, among others; see~\cite{automatedoutdoor,KR2020}.

The present survey paper is an expanded version of the keynote delivered at the ThEdu'22 workshop, August 2022, in Haifa (Israel).

\subsection{Computer Algebra Systems.}
\label{subsection CAS}
Decades ago, advanced specialized programs have been written, such as \href{https://www.gap-system.org/}{GAP - an acrostic for Groups, Algorithms, Programming} for Group Theory (its scope is much broader now), FeliX for Number Theory, \href{http://www2.macaulay2.com}{Macaulay2} for Algebraic Geometry, \href{https://cocoa.dima.unige.it}{CoCoA} for Computations in Commutative Algebra, etc. All these are examples of Computer Algebra Systems (CAS). General purpose CAS began also to integrate various mathematical fields, both symbolic and numerical algorithms; they are now powerful multi-domains assistants. The web offers also an interactive usage of platforms, where different CAS work in the background; this is the case with the education oriented platform \href{https://wims.univ-cotedazur.fr/}{WIMS}.

A core feature of a CAS is the ability to manipulate symbolic mathematical expressions, in a way similar to the traditional manual computations. Of course, the syntax of the commands can vary and some CAS may be thought as more user-friendly than  others. This can be an issue when working with students.

Among the needed features are a programming language (enabling the user to enter his/her own algorithms), an interpreter and a simplifier. A \textbf{simplify} command, with several options is important (pattern recognition may not be enough to lead to use an efficient algorithm). Note that for a given computation, the output may be quite different from what the user would have obtained by hand (if possible), and that two different packages may give different outputs. Simplifying helps the user to understand that different output formats may determine the same mathematical object. Some examples in Section~\ref{section examples} illustrate this fact.

Algorithms for Calculus, Linear Algebra, and other domains have to be used by students and teachers. Therefore, the CAS must include a large library of algorithms and special functions. We refer to a Wikipedia page\footnote{\url{https://en.wikipedia.org/wiki/Computer_algebra_system}} for a more exhaustive list of requested features, including some abilities of interest to the developers.
The first list on that page does not mention the plotting features: visualization, plotting, animations. They appear in a second list, not less important. We illustrate this in Section~\ref{section examples}.

The Wikipedia page adds that a CAS should also include a programming language, allowing users to implement their own algorithms, arbitrary-precision numeric operations, exact integer arithmetic and number theory functionality, editing of mathematical expressions in two-dimensional form, plotting graphs and parametric plots of functions in two and three dimensions, and animating them, drawing charts and diagrams,  APIs for linking it on an external program such as a database (we discuss later the need for a dialog between different systems), or using in a programming language to use the computer algebra system, string manipulation such as matching and searching, add-ons for use in applied mathematics such as physics, bioinformatics, computational chemistry and packages for physical computation, etc.
Some CAS include graphic production and editing such as computer-generated imagery and signal processing as image processing  and sound synthesis. DGS may fulfill a certain number of these requirements.

\begin{figure}[htb]
\centering
\mbox{
\epsfig{file=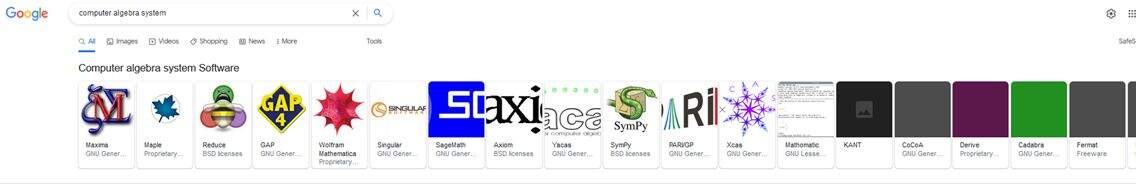,width=10cm}
     }
\caption{Google visual answer to a search for a CAS}
\label{fig google answer for CAS}
\end{figure}

As a visualization of the large variety of existing packages, Figure~\ref{fig google answer for CAS} shows the Google answer to a search for Computer Algebra System.

The danger is to use a CAS as a black box. We mean that the user writes a command, without any idea how it works. He/she enters the data, receives an output and relies on it. It is common sense that an educator will not go this way.

\subsection{Dynamic Geometry Systems.}
Later, other packages called Dynamic Geometry Systems (DGS) appeared: Geometric SketchPad, Cinderella, Cabri Geometer, GeoGebra, etc.. The needs for programming did not disappear, but a core feature is interactivity. The man-and-machine interaction is different from what is generally offered by a CAS.
A Dynamic Geometry System (DGS) provides ways of representing and manipulating geometric objects that are not possible with traditional paper-pencil work, using a compass and straightedge.
These various environments provide opportunities for students to explore geometric objects, to measure objects on the screen and outdoor. They can help students to develop different understandings of many properties and theorems. With the apparition of dragging and measuring offered by a DGS, teachers were sometimes reluctant to use them, as they feared that the role of the proof will fade. Mariotti~\cite{mariotti2001} reflected on these fears, but proposed a more positive view. Her reflection was based on work with the Cabri software.

A DGS  offers both button-driven commands and written commands, using a syntax similar to what exists in a CAS. A command has eventually a written version and a button, sometimes with slightly different affordances. For example, in GeoGebra's button for attaching a point to an object has a version as a written command working in cases the button does not. As a whole, this double feature buttons-written commands contribute to the software's user-friendliness.

The main dynamical features are the dragging points and slider bars. The basic objects are free points, and other constructions depend on these free points. Dragging a free point with the mouse induces automatically changes in every object depending on it. For example, the midpoint $I$ of a segment $AB$ depends on the endpoints $A$ and $B$. Dragging $A$ (or $B$, one at a time) with the mouse changes $I$. A slider bar corresponds to a real parameter; it allows to define and dynamically change objects depending on this parameter. For multi-parameter constructs, a slider is defined for each parameter (generally the software requests this automatically). Using sliders, animated graphs and constructs are available.  Dragging and sliders transformed Geometry into an experimental field. The main words in teaching Geometry are now exploration and experimentation. Changes using dragging and sliders, provide infinitely many examples of the situation under study. This is not a proof, but it provides conviction. Moreover, a conjecture can be enounced, and then has to be proven. The sequence exploration-conjecture-proof is ubiquitous when working with a DGS.

Animations are also offered by CAS, with an important difference. With a CAS, an animation is defined by a specific command with some options (size, number of frames per second, etc.). Then  the animation is run without an human intervention. The human can look at what happens, sometimes can intervene, but this is minor. An animation with a DGS can be almost automatic: the Animation On option available with a slider is an important feature, sometimes accompanied by the Trace On option for the object under study.  It can be mouse-driven, i.e. driven by the human. Examples are given in Section~\ref{section examples}. The interested reader can also refer to a dedicated Wikipedia page\footnote{\url{http://en.wikipedia.org/wiki/Interactive_geometry_software}}.

Mastering the moves with the hand makes the software some kind of a
prothesis for the human; Debray says that being able to build a
prothesis is what makes the man a human~\cite{debray}.  Maybe do we
deal now with an augmented human? We prefer to speak about a human,
totally human, with a strong digital literacy.

\subsection{Core issues with man-and-machine and machine-and-machine interaction}
\label{subsection machine interactions}

The availability of both CAS and DGS opened large new fields of study,
in research and in education. Their sets of respective affordances are
distinct, with a non empty intersection. Important results can be
derived by networking with the two kinds of
technology~\cite{roanes2003,DPKdialog}, one of them being more
powerful geometry and dynamics and the other having stronger algebraic
capabilities. After all, the study of a mathematical object is based
on different registers of representation (algebraic, numerical,
geometric, etc.)~\cite{duval2017}. For example, some packages are
devoted to 2D only or to 3D only. GeoGebra has both possibilities, the
2D window being fully synchronized with the $xy-$plane in the 3D
window.  Of course, they are also synchronized with the algebraic
window. Considering different representations of the same object,
i.e., multiple points of view, enable on the one hand to build bridges
between mathematical fields (v.i. Section~\ref{cognitive  neighborhoods}), and on the other hand, to derive new theorems or
new proofs of classical theorems (such as in~\cite{dpzrevival}). Until now, in order to have benefit of both environments,
CAS and DGS, we generally need to copy-paste from one system to the
other. The (lack of) communication between CAS and DGS has been
addressed for a long time~\cite{roanes2003,roanes2020,DPKdialog},
we wish that solutions will be found in the next future. Actually, one
of them was the embedding of a CAS into
GeoGebra~\cite{kovacs-parisse}, but as of today DGS and CAS are still
different. In the recent years, features of dynamic geometry entered
CAS and core CAS components have been implemented into
DGS~\cite{kovacs-parisse}. Work is still needed to have full
integration of all the possibilities.

Besides this machine-machine communication and collaboration, new
behaviours appeared in man-and-machine communication. This had also a
great influence on communication and collaboration between
humans. Shall we mention that the communication between humans and
between humans and machines developed new aspects during the Covid-19
crisis~\cite{open}?

The joint usage of the different technologies, including websurfing and data mining, enable huge advances in the various components of the so-called 4 C's of 21st Century Education: Communication, Collaboration, Critical Thinking and Creativity~\cite{4Cs}. The two first C's are generally considered between humans, but man-and-machine and machine-and-machine communication and collaboration are crucial. In~\cite{DPK-cassini-5C}, we propose even a 5th C, namely Curiosity, a must for human's proceeding further.  New didactic situations can be considered and new ways to deal with them. Mathematics Education has been transformed from the traditional Definition-Theorem-Proof scheme into an exploratory domain. The sequence can be now Exploration-Discovery-Conjecture-Proof (or disproof).

The 4C's are the conceptual basis of the STEAM Education approach. It
can be illustrated by examples and activities in Abstract Algebra,
Combinatorics together with Integrals, the study of Plane Algebraic
Curves and Algebraic Surfaces and also Mathematics and Arts. We refer
to the papers in~\cite{recioetalSI}.

\subsection{Numerical vs symbolic representations}
Between the symbolic data and the plot, a lot of numerical data is
computed, eventually displayed on demand. Numerical data induces also
the need to master the approximations~\cite{planetaryorbits}. This
numerical data is obtained, after choosing a mesh (in 3D) or a
partition of the interval of definition (in 2D), by
interpolation. This may lead to strange
plots~\cite{ZeDP2010,ZeDP2017}. This problem can be overcome either
with options for discontinuity or using a non standard mesh (Maple
proposes about 20 of non standard meshes for 3D plots).

\begin{figure}[htb]
\centering
\mbox{
\subfigure[]{\epsfig{file=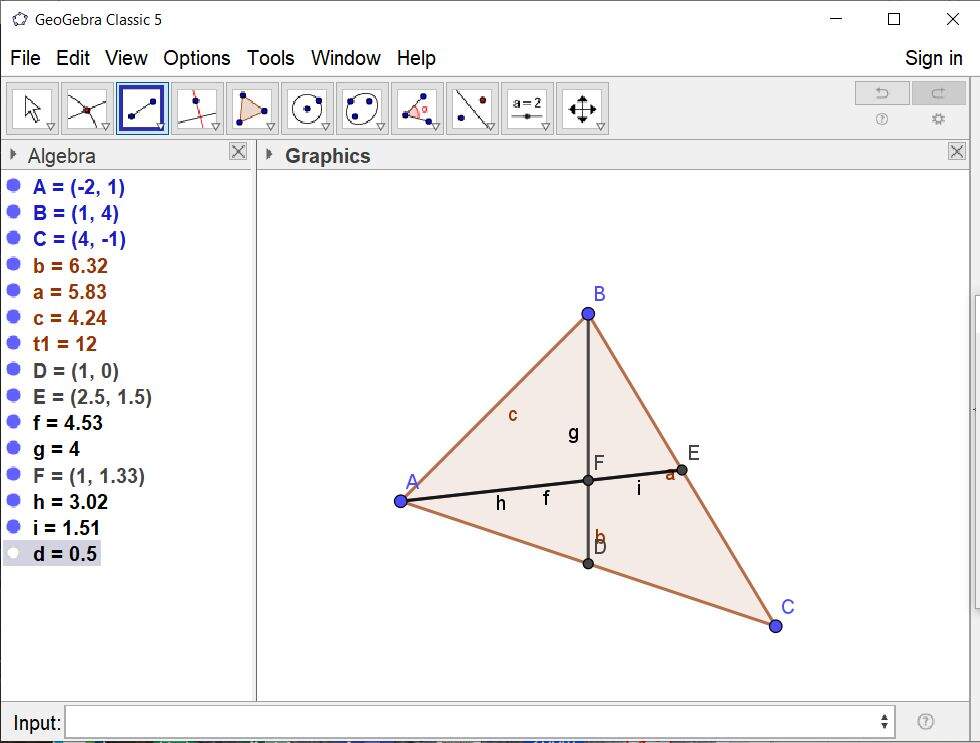,width=5cm}}
\quad
\subfigure[]{\epsfig{file=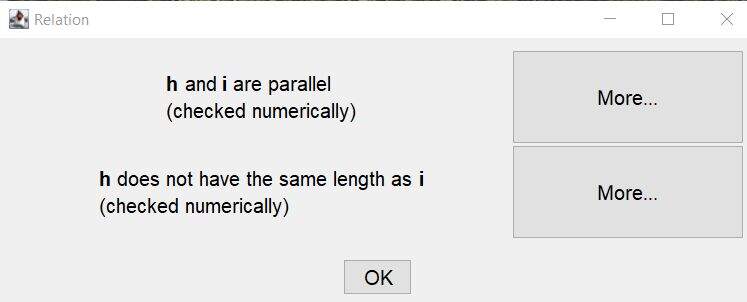,width=5cm}}
\quad
\subfigure[]{\epsfig{file=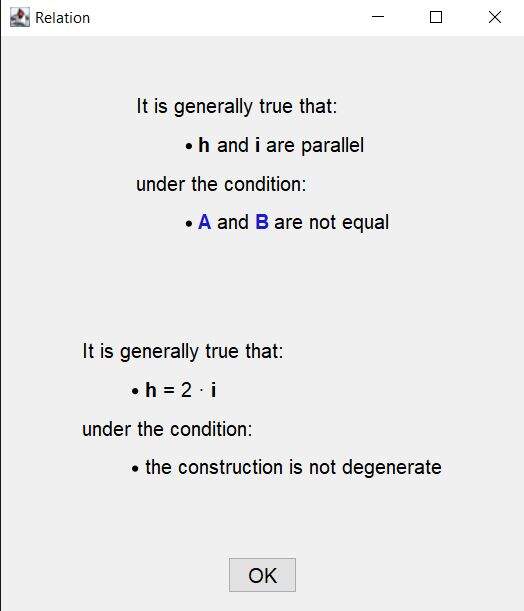,width=5cm}}
     }
\caption{Pop-ups offered by the Relation command}
\label{fig Relation}
\end{figure}

For other needs, the differences between numerical and symbolic algorithms led GeoGebra's developers to propose two kinds of output, for example, to the \textbf{Relation} command\footnote{We used here a beta-version of \emph{GeoGebra Discovery}, a package developed by Z. Kov\'acs, working on the basis of the regular GeoGebra version. It contains a more advanced version of the \textbf{Relation} command. The package is freely downloadable from \url{https://github.com/kovzol/geogebra-discovery} (check there for the last updated version).}. Figure~\ref{fig Relation} shows screenshots of the output.
The regular Relation command provides an answer based on numerical
data only. The ``More'' button runs other algorithms, that time
symbolic, the output is more precise.  It is an efficient tool in
Education, providing assistance to the teacher in developing new
activities.

\section{Cognitive neighborhood of a mathematical notion}
\label{cognitive neighborhoods}

Either with a single student or with a class, the teacher faces some endeavours
\begin{enumerate}
\item Stimulate  students' curiosity for  interlaced techniques, using more than one of the newly available technologies.
\item Make Mathematics  more attractive,  and show  it as  a living field of knowledge by discovering new tracks.
\item Discover links between apparently different fields. In a traditional curriculum, courses are  generally taught as separate topics, often without bridges between them. For example, some teachers are reluctant to introduce example for Calculus into a course in Linear Algebra, despite the natural structure of a vector space of the set of continuous (resp. differentiable) functions over an interval in $\mathbb{R}$.
\item Explore other mathematical objects ``looking  like'' (Consolidation?) the mathematical object of study. Such a task is generally built by the teacher, i.e. in this context, the teacher is active and creative.
\end{enumerate}

An advanced step in teaching integrals is provided by definite
integrals depending on one (ore more) parameter. This is a more
abstract situation than ordinary definite integrals, important in
applied situations~\cite{DPZ-soil}. The study of parametric definite
integrals provides intra-mathematics connections: integration,
combinatorics, applications to physical problems, and they may be
strengthened by topics in history of mathematics. The mathematical
notions, the bridges between them and the instrumented techniques,
accompanied by an appropriate technological discourse (a term coined
by Artigue~\cite{artigue}) build what we call a cognitive
neighborhood.  In the exploration of such a cognitive neighborhood,
teachers and students consider mathematical objects ``looking like''
the original mathematical object of study. Such a task is generally
built by the teacher, i.e. in this context, the teacher is active and
creative; the student reproduces the teacher's working steps. Such
activities incite the students to search for related material. The
exploration of the student's Zone of Proximal Development
(ZPD~\cite{vygotsky1978}) builds step by step an always larger
neighborhood.

Links to  neighboring mathematical topics, to ``real-world''  situations can be  discovered. Here the  student is more autonomous and can develop more initiative. Actually, both the educator and the student are creative.

\begin{figure}[htb]
\centering
\mbox{
\epsfig{file=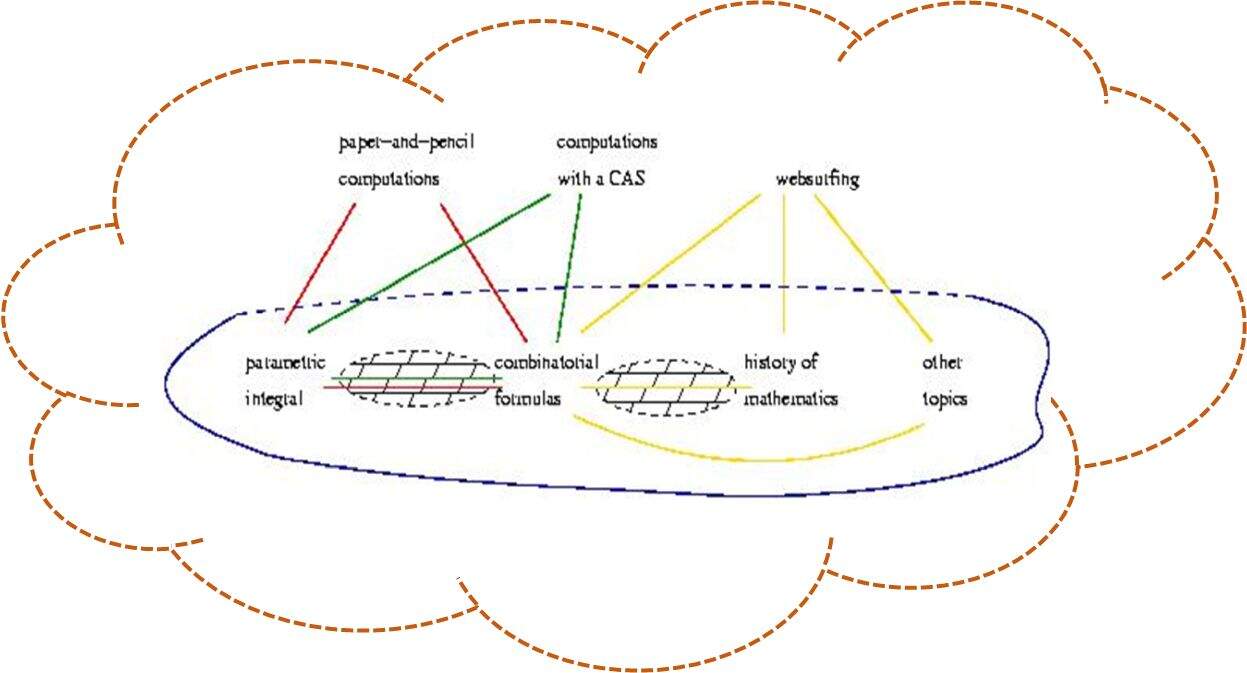,width=10cm}
     }
\caption{Visualization of a cognitive neighborhood}
\label{fig CognitiveNeighborhood}
\end{figure}

Figure~\ref{fig CognitiveNeighborhood} shows such a construct for
parametric integrals, as explored in~\cite{parametric2012,DPZe2017}.

\section{Next step: Automated Deduction in Geometry (ADG)}

This field of R \& D has seen tremendous developments during the last
decades, tens of papers and of conference presentations have been
devoted to the advances. Therefore we will mention only a small
sample. The interested reader can use the vast bibliography in the few
papers that we mention, and also have a look in the proceedings of
conferences such as ACA (Applications of Computer Algebra, especially
the special sessions on Education, but not only) and ADG
conference. In the ``far'' past, a short survey appeared
in~\cite{DPCASproofsthms2006}.

In~\cite{botana-abanades2014} Botana and Ab\'anades write the
following definition of ADG, ``as the study and development of computer
programs designed to prove geometry theorems''. Then they traced the
domain back to 1963 when Gelernter~\cite{gelernter} connected it to
Artificial Intelligence. They wrote: ``However, the real flourishing of
the field came in the early 1980's with the development by Wu of an
algebraic method based on Ritt's characteristic set for proving a
restricted set of geometry theorems~\cite{wu3}. Impressive results by
several authors using Wu's method~\cite{wu1,wu2,wu3} encouraged
researchers to consider other algebraic methods, among which those
based on Gr\"obner bases~\cite{buchberger} proved to be the most
relevant.''

Along the years, other researchers and developers joined the
domain. The definitions became more precise, in parallel with the
continuous software developments. In 2016, we could find the following
description~\cite{abanadesetal2016}: ``By automatic proving of
elementary geometry theorems, we refer to the theorem proving approach
via computational algebraic geometry methods, as initiated by Wu forty
years ago, and popularized by the book of Chou~\cite{chou}. Roughly
speaking, the idea is to provide algorithms, using computer algebra
methods, for confirming (or refuting) the truth of some given
geometric statement. More precisely, the goal is to decide whether a
given statement is generally true or not, i.e. true except for some
degenerate cases, to be described by the algorithm.''

The developments of different kinds of software, in the present case
DGS and Automatic Provers enabled automated reasoning and dynamic
geometry together. The collaboration of the different software
revealed very successful. In 2020, Kov\'acs and Recio
write~\cite{KR2020}: ``Along the last half century, automated deduction
in elementary geometry has increasingly become one of the most
successful achievements in the field of automated reasoning. Along
these decades various methods and techniques have been studied and
developed for automated proving and discovering of elementary geometry
statements. On the other hand, dynamic geometry software systems have
emerged, such as Cabri Geometry, C.a.R., Cinderella, DrGeo, GeoGebra,
The Geometer's Sketchpad, Geometry Expert, Geometry Expressions or
Kig, with an ever-increasing presence in mathematics education. Some
of them possess a large number of users (over thirty million) all
around the world. The merging of these two tools (automatic proving
and dynamic geometry) is, thus, a very natural, challenging and
promising issue, currently involving logic, symbolic computation,
software development, algebraic geometry and mathematics education
experts from all over the world.''

A central topic in Geometry is the determination of geometric
loci. Numerous results have been obtained, some of them are presented
in~\cite{botanaloci,locusblazek}.  We explore some loci in
subsection~\ref{subsection loci}. Commands for automated exploration
and discovery of geometric loci have been developed and implemented in
GeoGebra. We illustrate them in subsection~\ref{subsection loci}.

In subsection~\ref{subsection envelopes}, we address another topic,
namely envelopes of parametric families of plane curves. Despite
Thom's complained~\cite{thom} (back in 1962) that it disappeared from
the syllabus, research did not stop and numerous papers have been
devoted to envelopes, software developments and Mathematics
Education. A small sample
is~\cite{botana-valcarce2,botanarecioJSSC,dps1996,botanarecioAMAI,safety,DPK-cassini-5C}. For
practical situations, approximate methods had also to be developed;
for example, see~\cite{schultz-juttler,pottmann-peternell}. An
automated command has been implemented in GeoGebra for the
determination of envelopes of families of plane curves. Examples are
given in subsection~\ref{subsection envelopes}. The first computations
consist in solving a non-linear system of equations, which is feasible
using a CAS. The output consists in a list of parametric presentations
of curves, and may need to be simplified, whence the importance of the
simplify command. Examining the equations together with plotting them
reveal that the components are complementary. This is to be proven,
generally with algebraic tools. First the parametric presentations
have to be translated into polynomials, then Gr\"obner packages can be
applied. These packages may not be implemented in the DGS, and work
has to be transferred to another CAS, generally using copy-paste.  The
obtained polynomials generate an ideal in a polynomial ring. By
elimination of the parameter, an implicit equation can be derived for
the envelope. This is performed with a package for Gr\"obner bases
computations; see~\cite{cox,montes} for the theory of Gr\"obner bases
and~\cite{sendra} for rational curves. Using once again the CAS, the
obtained polynomial can be factored, revealing whether the envelope is
reducible or not. Of course, the CAS offers plotting features, and
their output will give a confirmation to what has already been
obtained with the DGS.

Different definitions of an envelope are given
in~\cite{bruceandgiblin}; we illustrate them with automated methods in
subsection~\ref{subsection envelopes}. Offsets are worked out as loci,
the difference between them and envelopes are made clear by using a
DGS, as in~\cite{talbot}.

The third topic that we will address later is the determination of
isoptics curves of plane curves.  It has been studied for
years~\cite{cieslak-miernowski-mozgawa1991,miernowski-mozgawa1997,szalkowski}. Working
in a technology-rich environment enabled new
results~\cite{DMZ-bisopticsofellipses,isopticsfermat}. More
recently inner isoptics have been studied~\cite{dp-mozgawainner}. Some
dynamics have been introduced into the study of isoptic
in~\cite{DPK-dynamicalcoloring}. As of today, no specific command has
been implemented for the discovery of isoptics; nevertheless the
existing features allow exploration, computations and confirmation of
the results. This is the topic of subsection~\ref{subsection
  isoptics}.

To close this section, we wish to mention the book by Pavel
Pech~\cite{pechbook}. A chapter is devoted to automatic theorem
proving and automatic discovery. The subsequent chapters deal with
classical theorems in Geometry by means of polynomials with CoCoA, a
freely downloadable software for computing in polynomial rings;
see~\cite{cocoalib}.

\section{Other examples}
\label{section examples}
We present now a couple for mathematical activities for which technology may have a crucial role. For exploration, discovery and proofs, we use here GeoGebra and Maple.

\subsection{Experiment the continuity of a function}
\begin{figure}[htb]
\centering
\mbox{
       \subfigure[continuity]{\epsfig{file=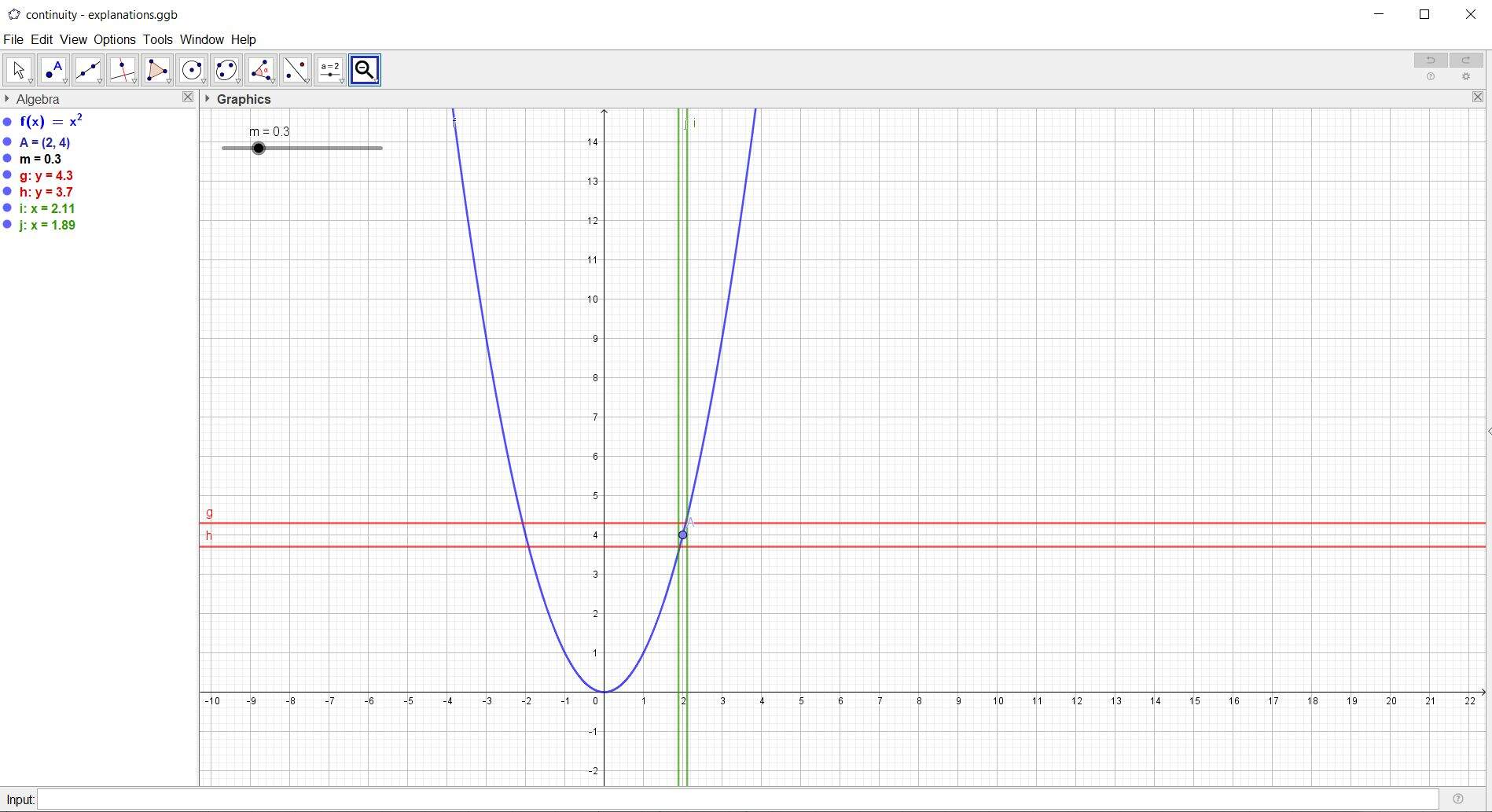, width=5.5cm}}
       \quad
       \subfigure[discontinuity]{\epsfig{file=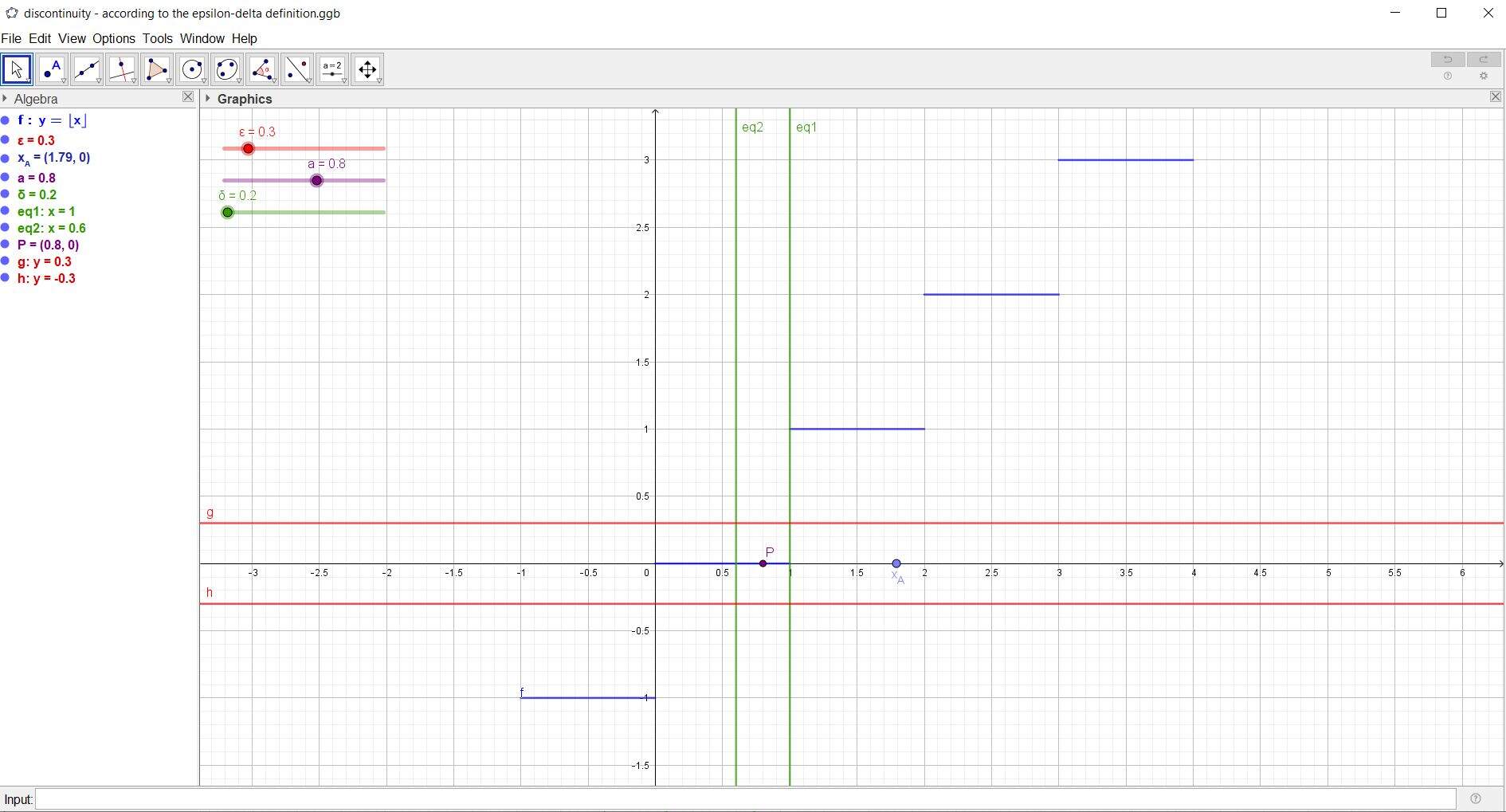, width=5.5cm}}
     }
     \caption{Exploration of the continuity of a function according to the  $\epsilon - \delta-$ definition}
     \label{fig applet continuity}
\end{figure}

\begin{figure}[htb]
\centering
\mbox{
\epsfig{file=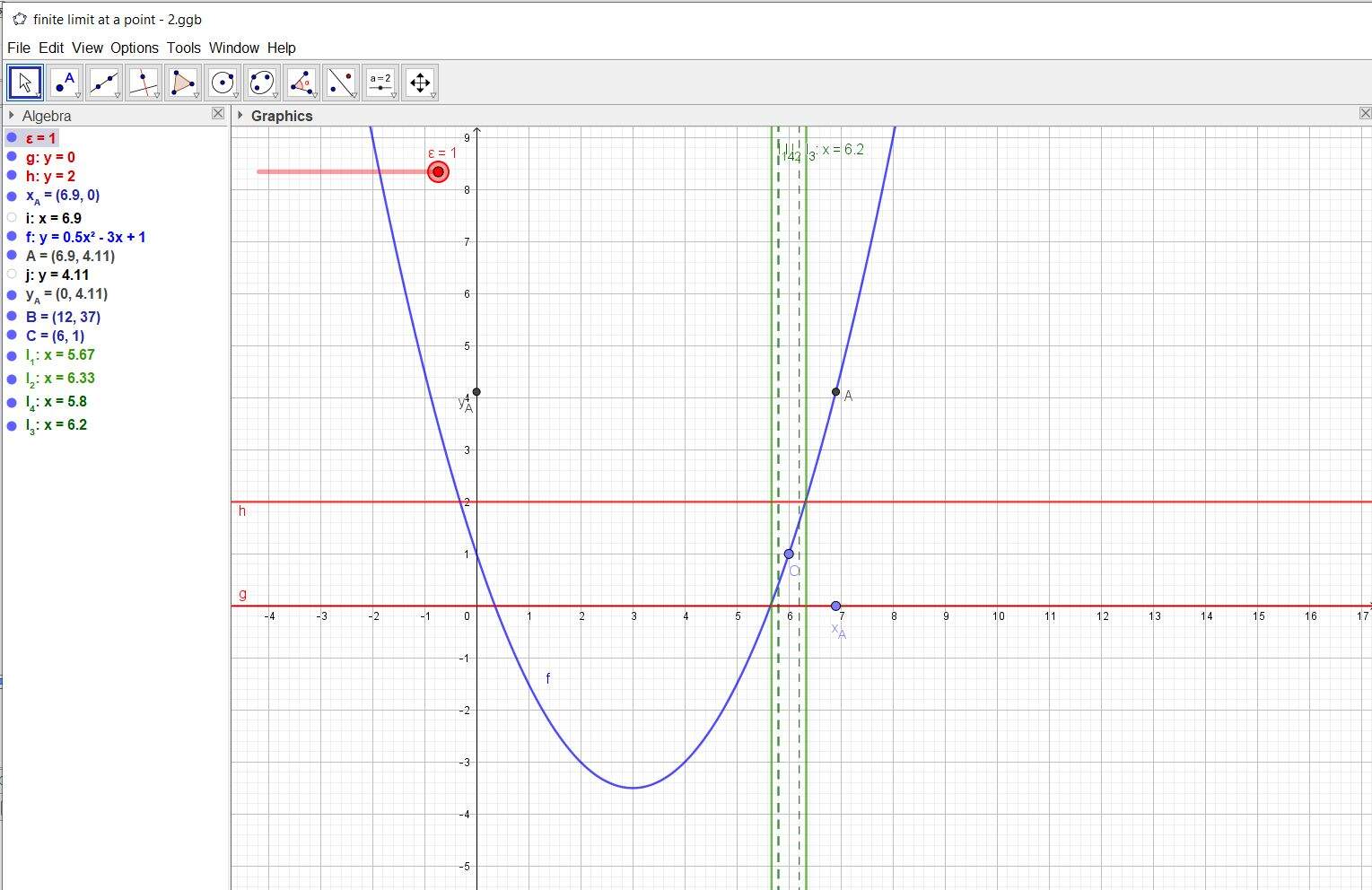,width=5.5cm}
     }
\caption{A more advanced step for illustration of continuity}
\label{fig applet continuity advanced}
\end{figure}

Mastering the  $\epsilon - \delta-$ definition is often difficult for a student. Exploration with technology makes things easier to understand. In order to make things easier, more intuitive, we developed an applet with GeoGebra\footnote{Available at \url{https://www.geogebra.org/m/bkj8872v}}. Figure~\ref{fig applet continuity} shows two screenshots. The sliders enable to change independently the value of $\epsilon$ and of $\delta$, which is useful for a first exploration. The student notes that $\epsilon$ determines a horizontal band and $\delta$ a vertical one, whose intersection is  a rectangle centered at the point $(x_0,f(x_0)$, the question being whether the graph of the function passes fully in the rectangle. Discontinuity appears when $\epsilon$ is quite small and no value of $\delta$ enables this.
Later, the teacher may propose another applet\footnote{\url{https://www.geogebra.org/m/zsb6p2vd}}, where the dependence of $\delta$ on $\epsilon$ is utilized, as shown in Figure~\ref{fig applet continuity advanced} and only one slider appears. This applet requires preliminary human work to discover a formula for the dependence of $\delta$ on $\epsilon$. Therefore it is not suitable for exploration, more suitable for illustration, and even be useful for consolidation of the new knowledge. Of course, the best situation is when the student is able to develop such an applet, and does not rely on the teacher for this. A student may develop suitable skills for this. If not the teacher may offer some scaffolding (\cite{noss-hoyles}, pages 105-107).

\subsection{Exploration of a geometric locus}
\label{subsection loci}
GeoGebra has several versions of an automated command for the determination of a geometric locus. According to the situation, whether the tracer is geometrically related to the mover or the tracer has to fulfill a given condition (a Boolean expression), the output has different forms: it can be a plot of a curve or a plot together with an implicit equation. Recent developments provide plots of regions in the plane\footnote{Derive did it, representing inequalities in 2 variables graphically.}.

\subsubsection{Classical and less classical constructions.}
Take two points $A$ and $B$ in the plane. The geometric locus of points $M$ such that the ratio $AM/BM$ is equal to $k$ is the perpendicular bisector of the segment $AB$ if $k=1$, and a circle otherwise. This can be explored using the Locus command and a slider for letting the ratio vary (see \url{https://www.geogebra.org/m/vqqkd57t}). For each value of the parameter an implicit equation is provided.

\begin{figure}[htb]
\centering
\mbox{
\subfigure[]{\epsfig{file=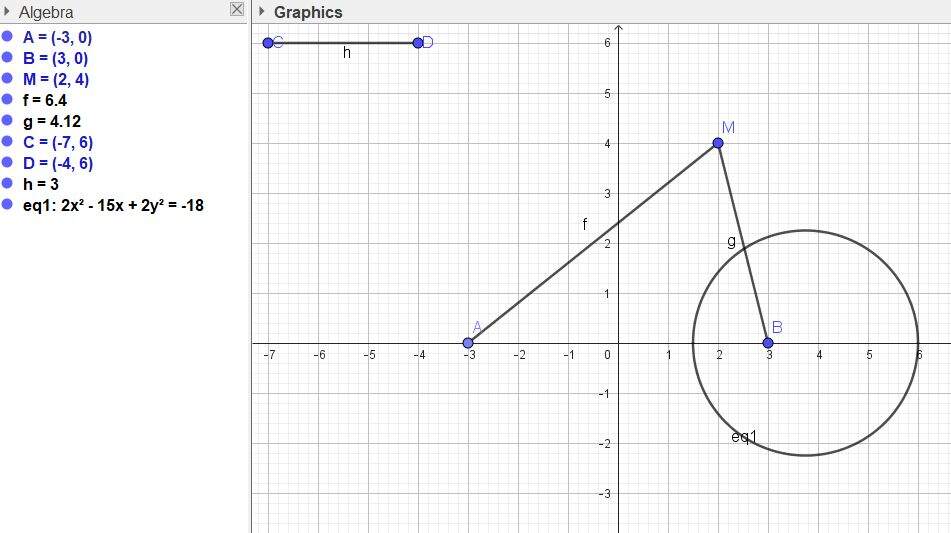,width=5cm}}
\qquad \qquad
\subfigure[]{\epsfig{file=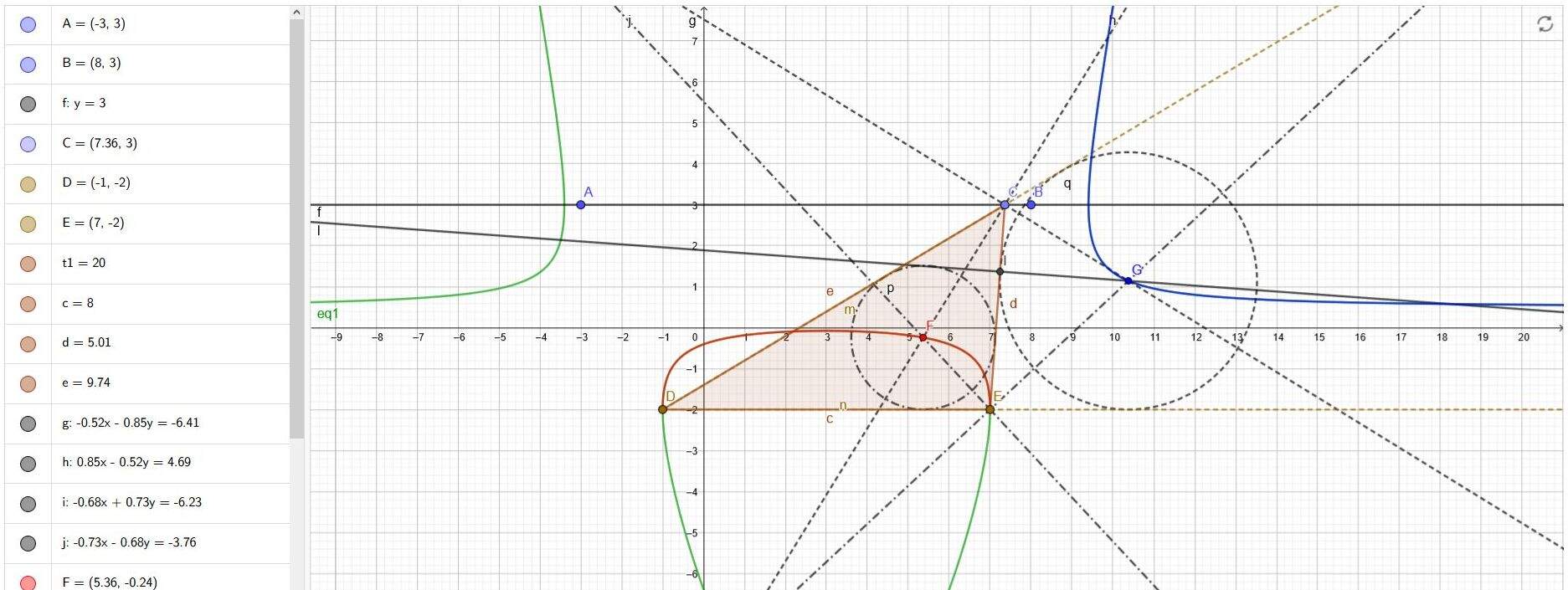,width=7cm}}
}
\caption{Geometric loci}
\label{fig loci - geometric constructions}
\end{figure}

The applet  \url{https://www.geogebra.org/m/v2kdhrus} shows the construction of the loci of the incenter and excenters of a triangle, where two vertices are fixed and the third one moves along a line; see Figure~\ref{fig loci - geometric constructions}(b). Here the version of the command is Locus($<$Point Creating Locus line $>,<$Point$>$) and the output is a plot of the locus, but no equation is found.

The same kind of construction can be applied looking for the geometric locus of points $M$ such that the angle $\angle AMB=\theta$ fro given $\theta$. The output is a circle. In this case, the teacher has to draw the students' attention that, actually, the geometric locus is the circle without the two points $A$ and $B$, as if $M \in \{ A,B \}$, there is no angle at all. A good opportunity to educate to Critical Thinking~\cite{4Cs}.

\subsubsection{Beyond the directrix of a parabola.}
\begin{figure}[htb]
\centering
\mbox{
\subfigure[]{\epsfig{file=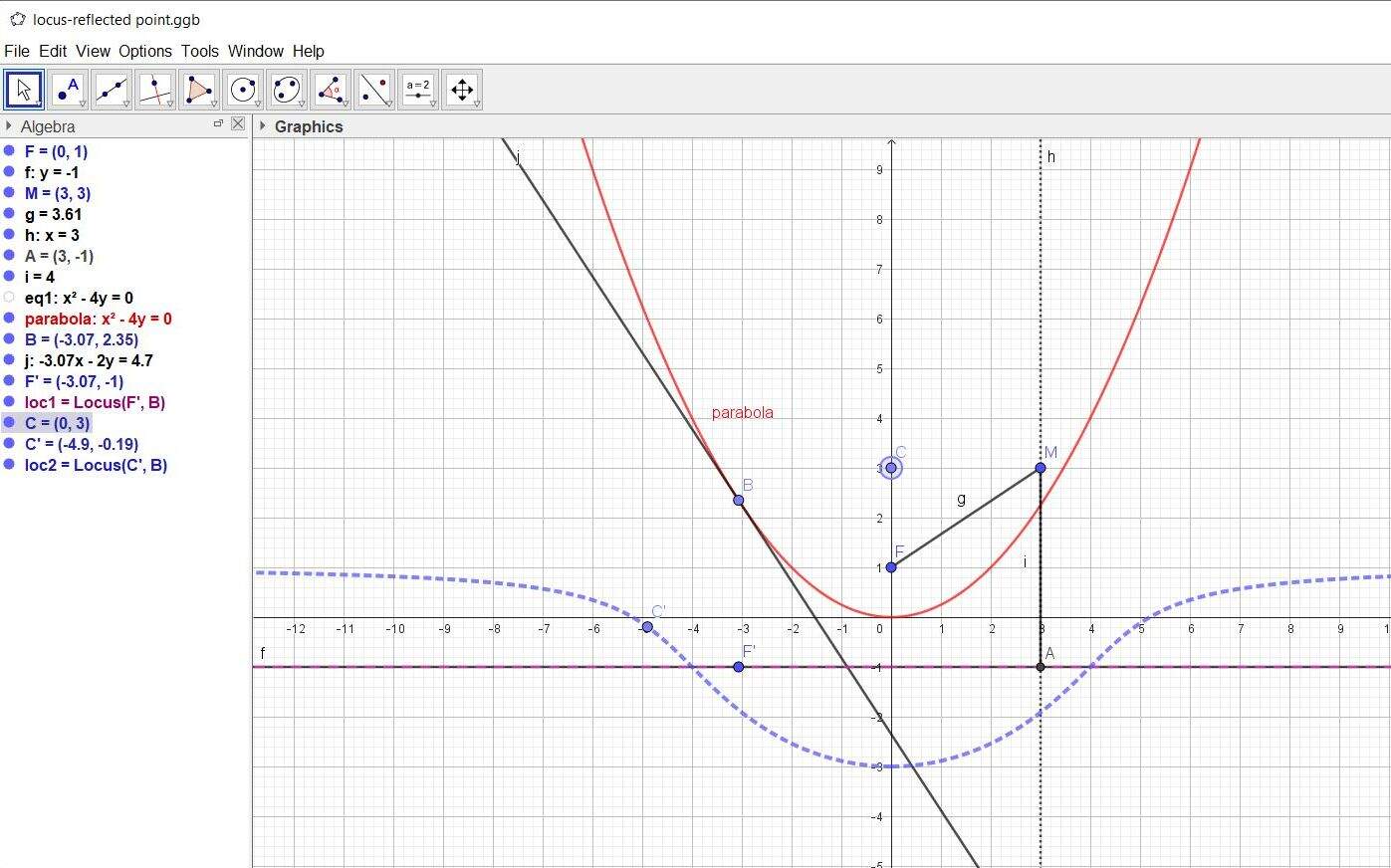,width=5cm}}
\qquad \qquad
\subfigure[]{\epsfig{file=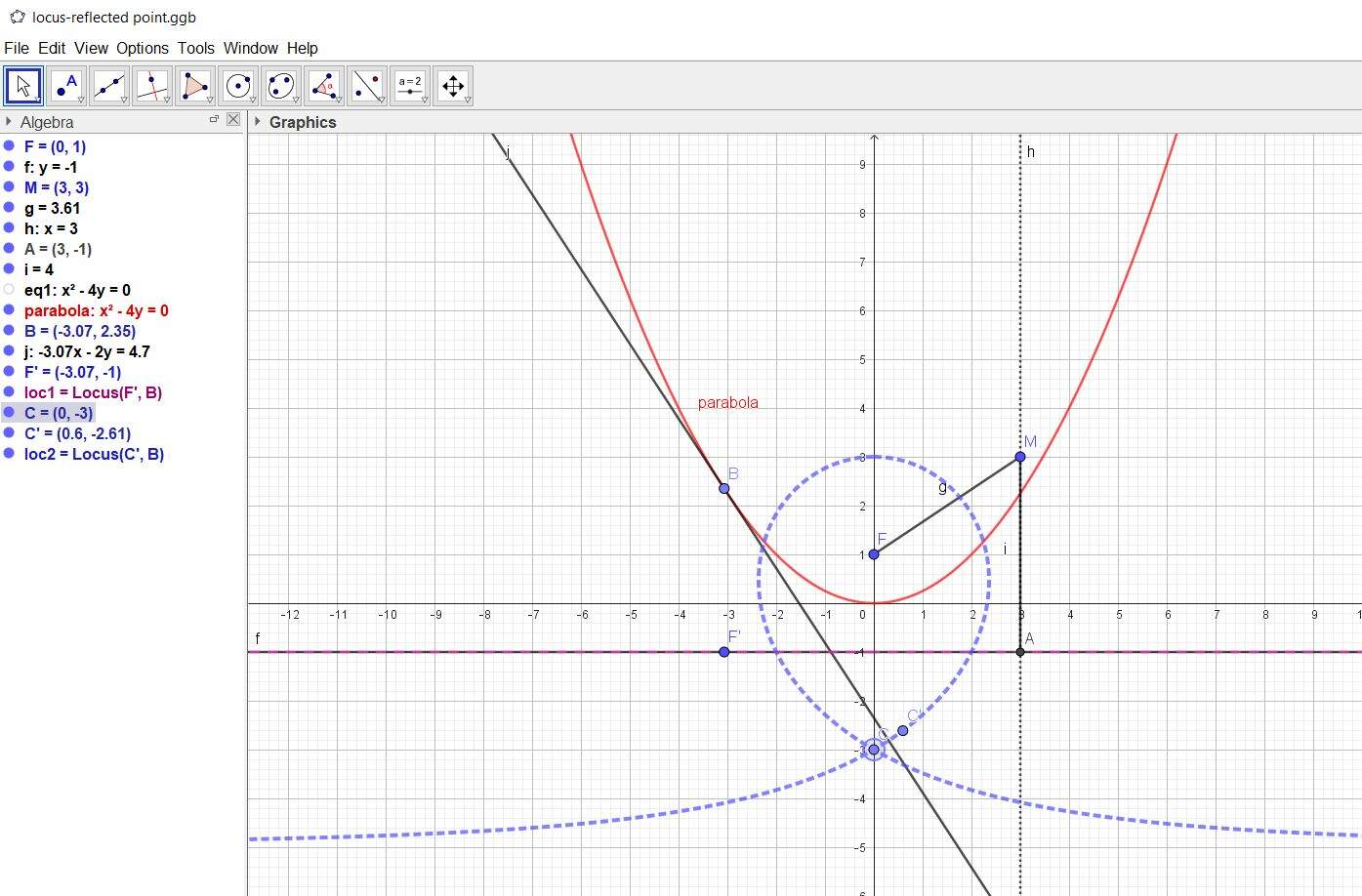,width=5cm}}
}
\caption{The locus of reflected parabola $y-$axis point about tangents}
\label{fig locus of reflected parabola y-axis points}
\end{figure}

Let $\mathcal{P}$ be the parabola whose equation is $y=x^2$. Denote by  $B$  a point on $\mathcal{P}$ and by $C$ a point on the $y-$axis of $\mathcal{P}$. Let $C'$ be the image of $C$ by the reflection about the tangent to $\mathcal{P}$ at $B$. Two examples are on display in Figure~\ref{fig locus of reflected parabola y-axis points}.
The Locus($<$Point Creating the Line$>,<$Point$>$) creates the dotted plot without providing an equation. A database of plane curves may be searched in order to determine which curve has been obtained, and to try to find a suitable equation for it. Students observed that the limiting case between smooth curves and curves with a double point is when the point $C$ is the focus of the parabola $\mathcal{P}$, in which case the geometric locus of $C'$ is the directrix of the parabola.

 \subsubsection{A crosscurve.}
\begin{figure}[htb]
\centering
\mbox{
\epsfig{file=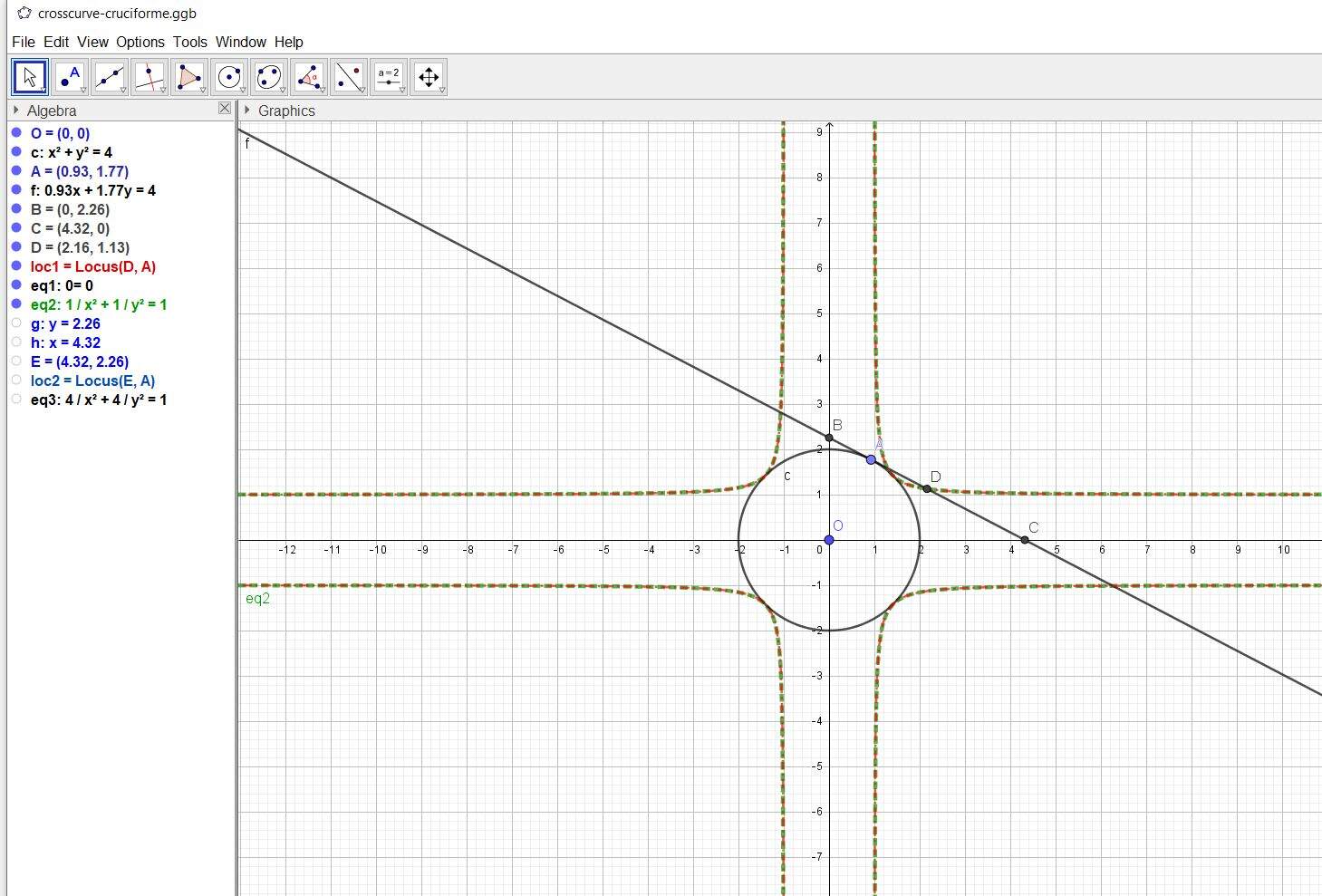,width=5.5cm}
     }
\caption{Automated exploration of a geometric locus}
\label{fig crosscurve polynomial}
\end{figure}

Let be given a circle centered at the origin. We explore the geometric locus of the midpoint of the intercepts with the axes of the tangents to the circles. An automated command is available in GeoGebra, the Locus command. Several version are proposed, and the student has to choose among them.  A snapshot of  a session\footnote{\url{https://www.geogebra.org/m/bmdfd3hm}} is displayed in Figure~\ref{fig crosscurve polynomial}:  a purely geometric construct is performed, enabling the usage of the Locus ($<$Point$>,<$Point$>$) command. Here the tracer is the midpoint and the mover is the point on the circle (for details on the terminology see~\cite {ecosystem}).
Identification of the obtained curve may not be easy. In our example, the circle has radius 2. The geometric locus has equation
\begin{equation}
\label{eq crosscurve}
\frac {1}{x^2}+\frac{1}{y^2}=1
\end{equation}
This equation has been derived by either by hand or with a CAS. In reverse direction, it could have been discovered in a catalogue of classical plane curves\footnote{Such as the Mathcurve website. The curve appears there at \url{https://mathcurve.com/courbes2d.gb/cruciforme/cruciforme.shtml}}, but this an unilluminating search. Maybe a websearch based on a picture will help, we did not try it.

\begin{remark}
 Equation (\ref{eq crosscurve}) can be transformed into a polynomial equation, namely
 \begin{equation}
\label{eq crosscurve}
x^2+y^2=x^2y^2.
\end{equation}
This may be an indication that packages for polynomial computations may be useful for the current exploration. Here too, Critical Thinking has to be applied as the domains of validity  are different.
\end{remark}

\subsubsection{A sextric.}
\begin{figure}[htb]
\centering
\mbox{
       \subfigure[$r=1/3$]{\epsfig{file=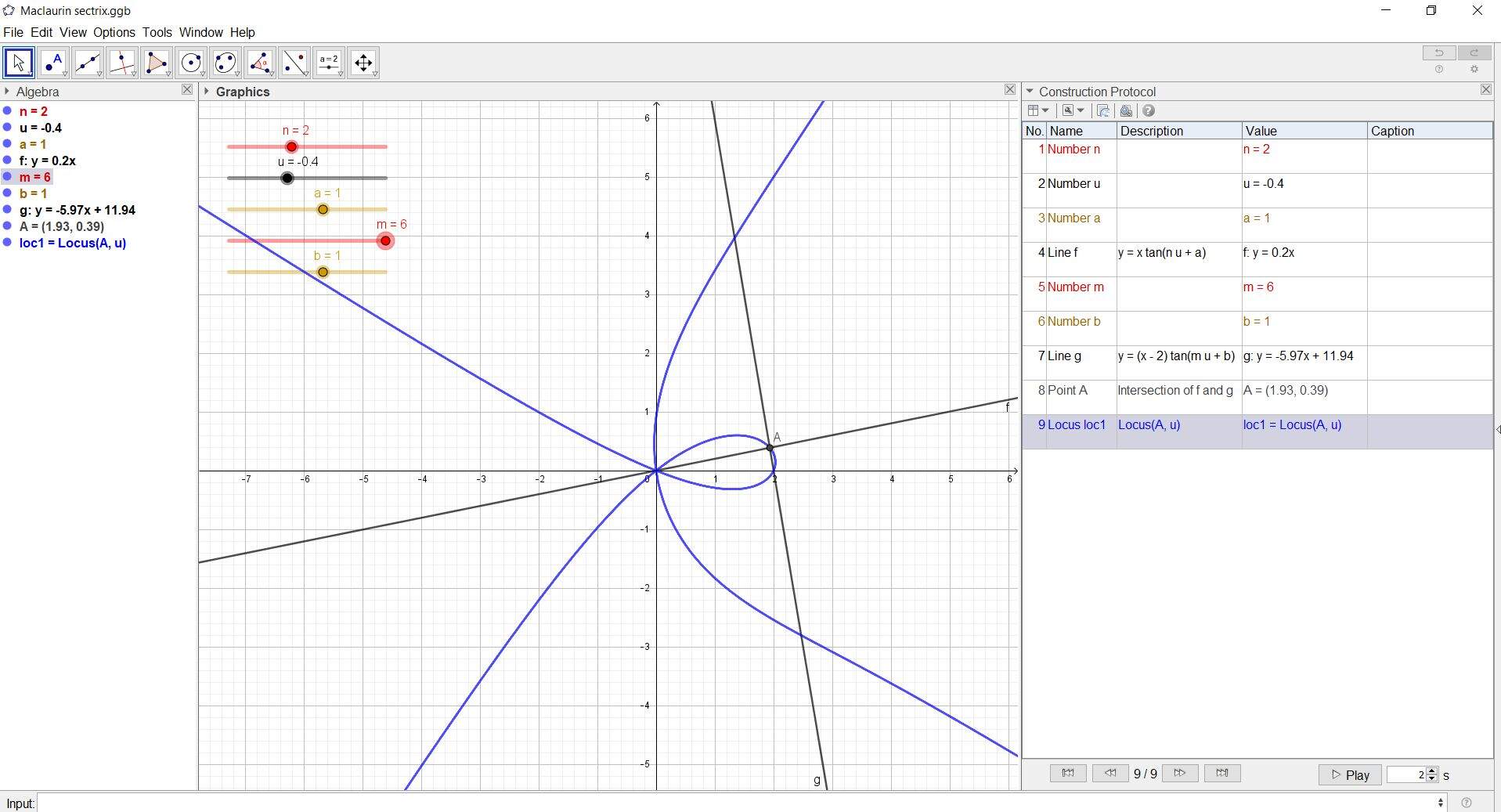, width=5cm}}
       \quad
       \subfigure[$r=-1$]{\epsfig{file=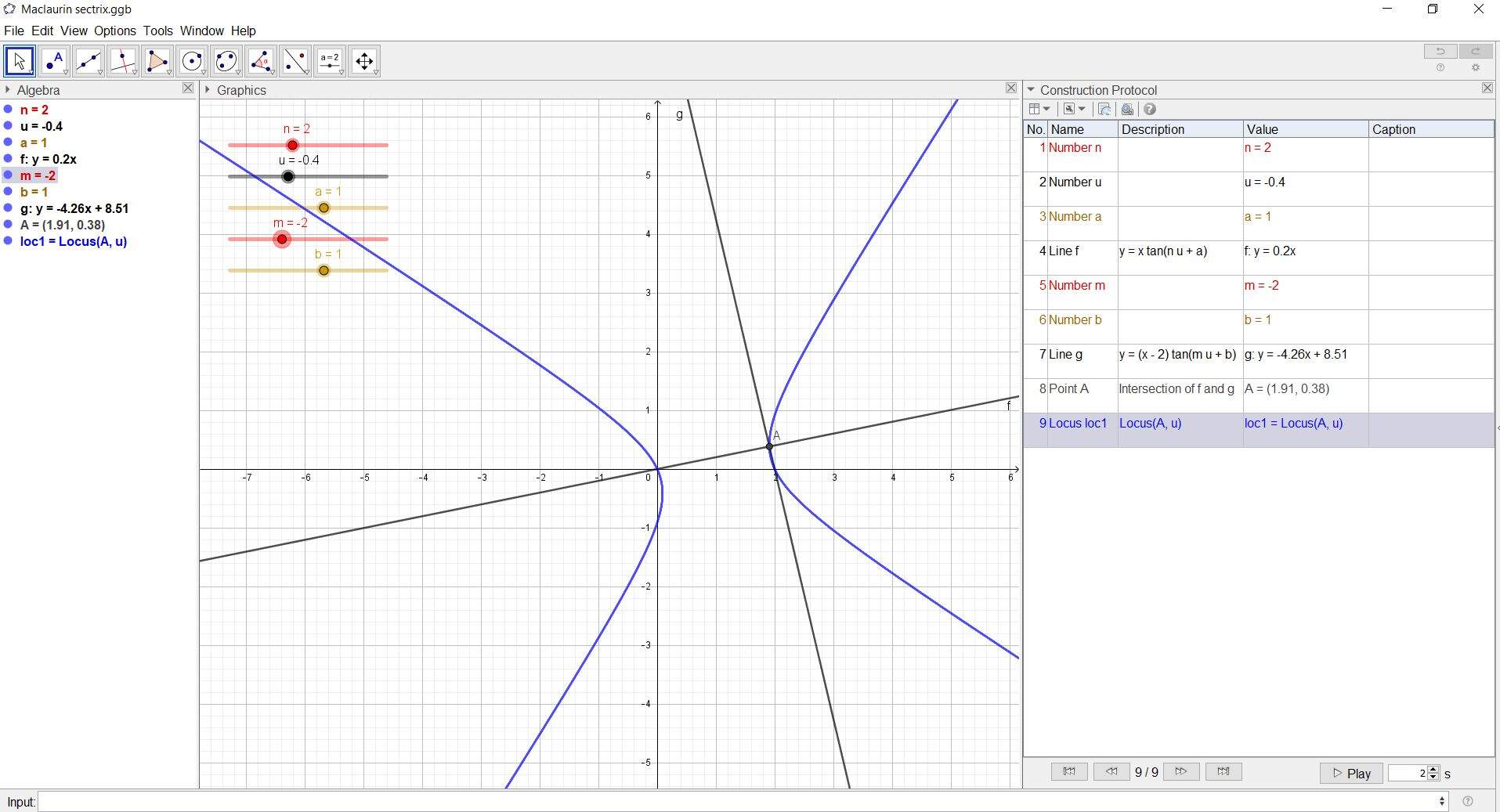, width=5cm}}
       \quad
       \subfigure[$r=2$]{\epsfig{file=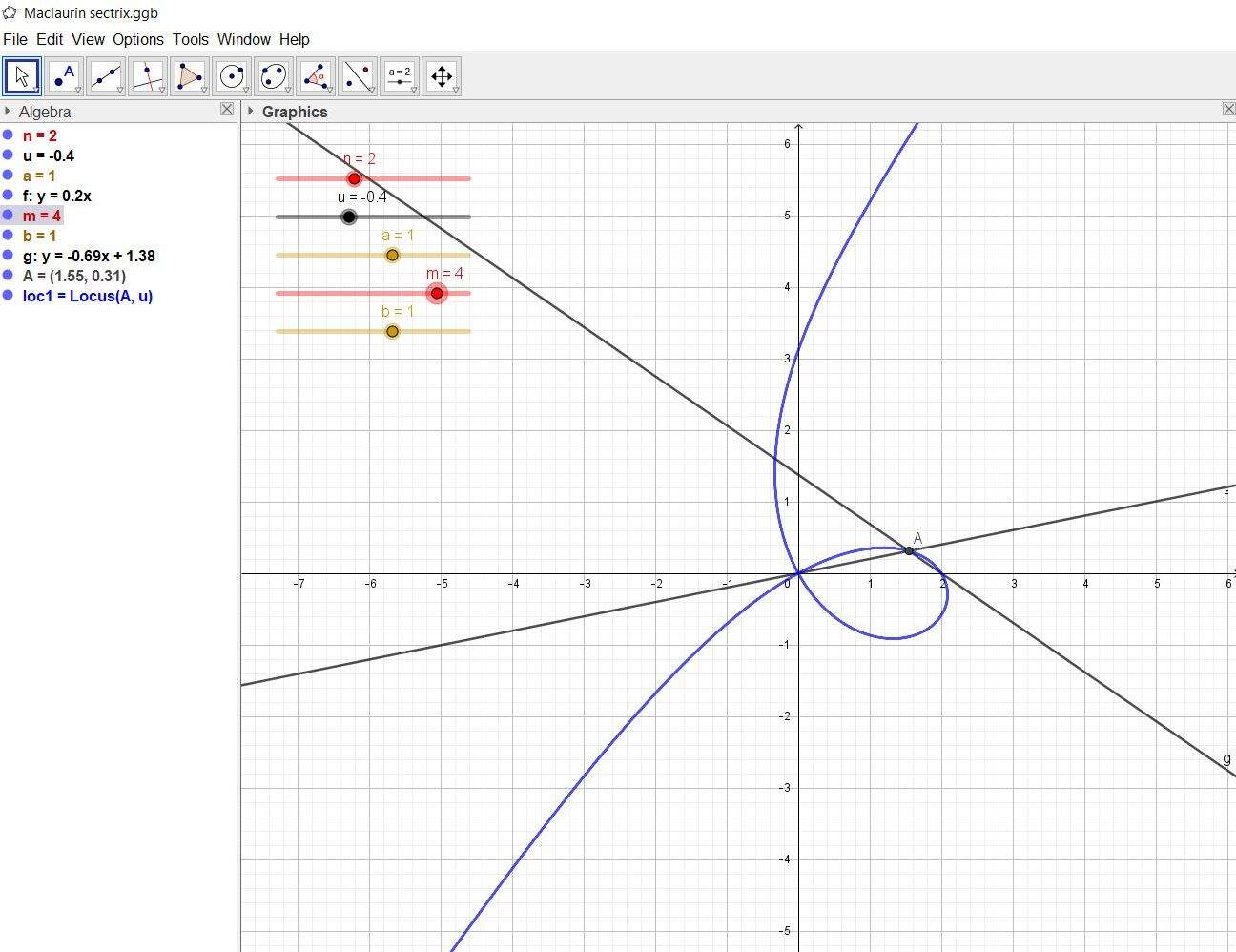, width=5cm}}
       }
     \caption{sextrices of MacLaurin}
     \label{fig sextric of MacLaurin}
\end{figure}

A sextric is a curve of degree 6. A \emph{sextric of Maclaurin} is the geometric locus of the point of intersection of two lines, which are each revolving at constant rates about different points called poles. Figure~\ref{fig sextric of MacLaurin} shows three examples, according to the ratio $r$ of angular velocities being equal to 1/3, -1 or 2 (or, equivalently, to 1/2). In the lst case, the obtained curve is a strophoid. Note the different numbers of components. These are screenshots from a GeoGebra applet\footnote{\url{https://www.geogebra.org/f/kwzvga3rpn}}; exploration is made possible by the definition of several sliders.

Catalogues of plane curves show them as separate objects. The slider provides continuous changes in the values of the parameter, leading to discover that several seemingly different curves belong actually to a single family. The apparition of a straight component requires more theoretical developments.

\subsection{Envelopes}
\label{subsection envelopes}
There exist 4 different definitions of an envelope;
see~\cite{bruceandgiblin} (Chap. 5). Kock~\cite{kock} gives three of
them, calling them respectively synthetic, impredicative and
analytic. We will explore the examples, according to the different
definitions, and see how automated methods are applicable. More
advanced examples may be found in~\cite{dpzrevival,safety,talbot};
dedicated databases provide also numerous examples.

\begin{definition}
\label{def envelope synthetic}
\textbf{[Synthetic]} Let $\mathcal{C}_u$ be a family of real plane curves dependent on a real parameter $u$.
 The envelope $\mathcal{E}$ is the union of the characteristic points $M_u$, where the characteristic point $M_u$ is the limit point of intersections $\mathcal{C}_u \cap \mathcal{C}_{u+\epsilon}$ as $h \rightarrow 0$. In other words, the envelope is the set of limit points of intersections of nearby curves $\mathcal{C}_u$,
\end{definition}

\begin{definition}
\label{def envelope impredicative}
\textbf{[Impredicative]} The envelope $\mathcal{E}$ is a curve such that at each of its points, it is tangent to a unique curve from the given family. The locus of points where $\mathcal{E}$ touches $\mathcal{C}_u$ is called the $\mathcal{E}-$characteristic point $M_u$.
\end{definition}

\begin{definition}
\label{def envelope analytic}
\textbf{[Analytic]} Suppose that the family of curves $\mathcal{C}_u$ is given by an equation $F(x,y,u)=0$ (where $u$ is a real parameter and $F$ is differentiable with respect to $u$), then an envelope $\mathcal{E}$ is determined by the solution of the system of equations:
\begin{equation}
\label{system def envelope}
\begin{cases}
F(x,y,u)=0 \\
\frac{\partial F}{\partial u} F(x,y,u)=0
\end{cases},
\end{equation}
\end{definition}
In Definition~\ref{def envelope analytic}, the envelope is described
as the projection onto the $(x, y)$-plane of the points, in the
3-dimensional $(x, y, u)$-space, belonging to the surface with
equation $F(x, y, u) = 0$ and having tangent plane parallel to the
$u$-axis (or being singular points and, thus, not having tangent
plane, properly speaking). See~\cite{bruceandgiblin}, p.102. Note that
the analytic definition~\ref{def envelope analytic} is the only one
given by Berger~\cite{berger}(sections 9.6.7 and 14.6.1) and by
Rovenski~\cite{rovenski}. This last book gives details on how to work
out envelopes using Maple.

\begin{remark}
We chose to denote the parameter by $u$, as in GeoGebra $t$ has a special role.
\end{remark}

\subsubsection{An envelope of a family of lines}
We consider a family of lines given by the equation $F(x,y,u)=0$,
where $F(x,y,u)=x+uy+u^2$. Such a family has been studied with Derive
in~\cite{DMZ-bisopticsofellipses}; we study it here in order to show
the usage of automated commands, if possible. Figure~\ref{fig
  lines-CAS}(a) shows a first exploration, using the slider and Trace
On\footnote{\url{https://www.geogebra.org/f/qt5pdbah42}. Note that the
  increment for the slider $a$ has to be put to 0.05 in order to have
  an accurate plot.}. The speed of the mouse on the slider conditions
the density of the lines in the output. Not plotting too many lines
retains the visual impression of a family of lines, and not a fully
colored area.

\begin{figure}[htb]
\centering
\mbox{
\subfigure[With Trace On]{\epsfig{file=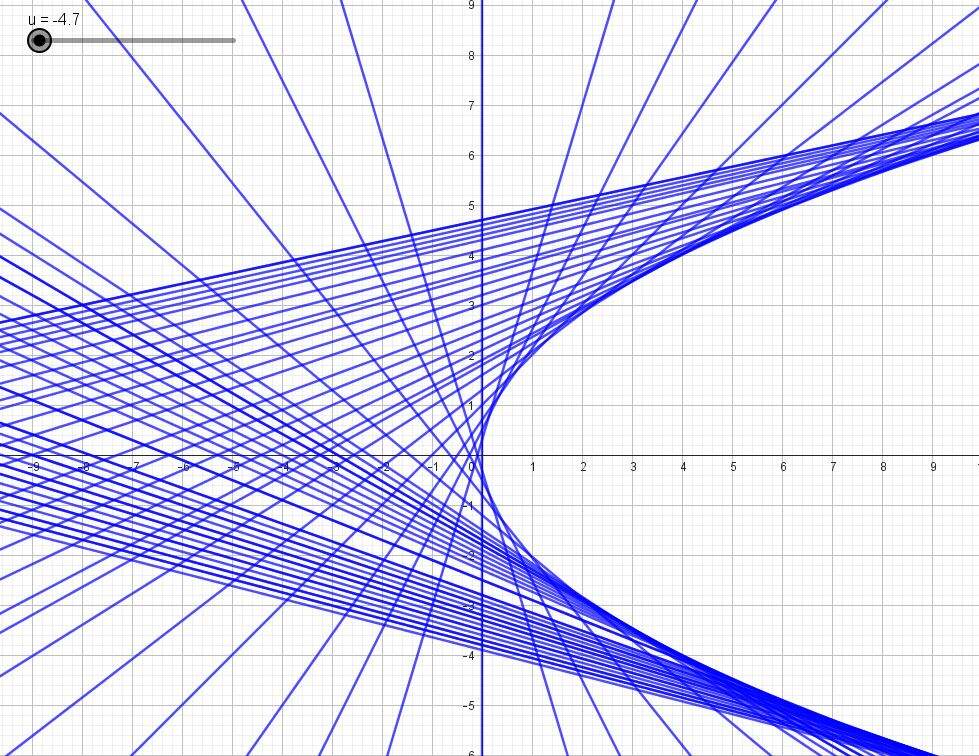,width=4.5cm}}
\quad
\subfigure[Two neighboring lines]{\epsfig{file=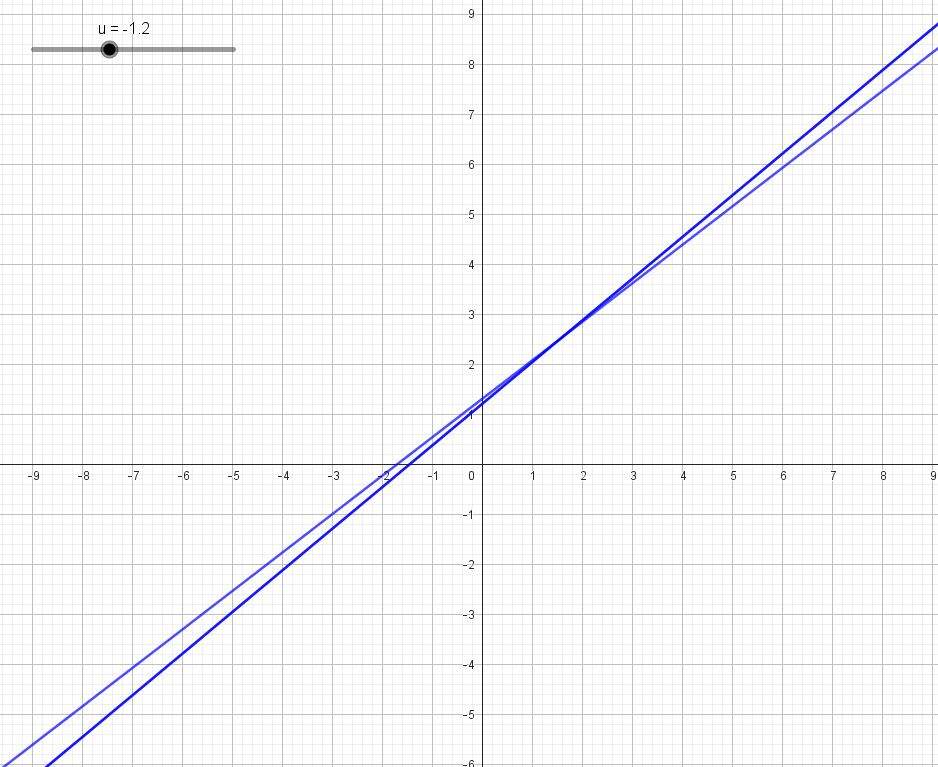,width=4.5cm}}
\quad
\subfigure[With the envelope]{\epsfig{file=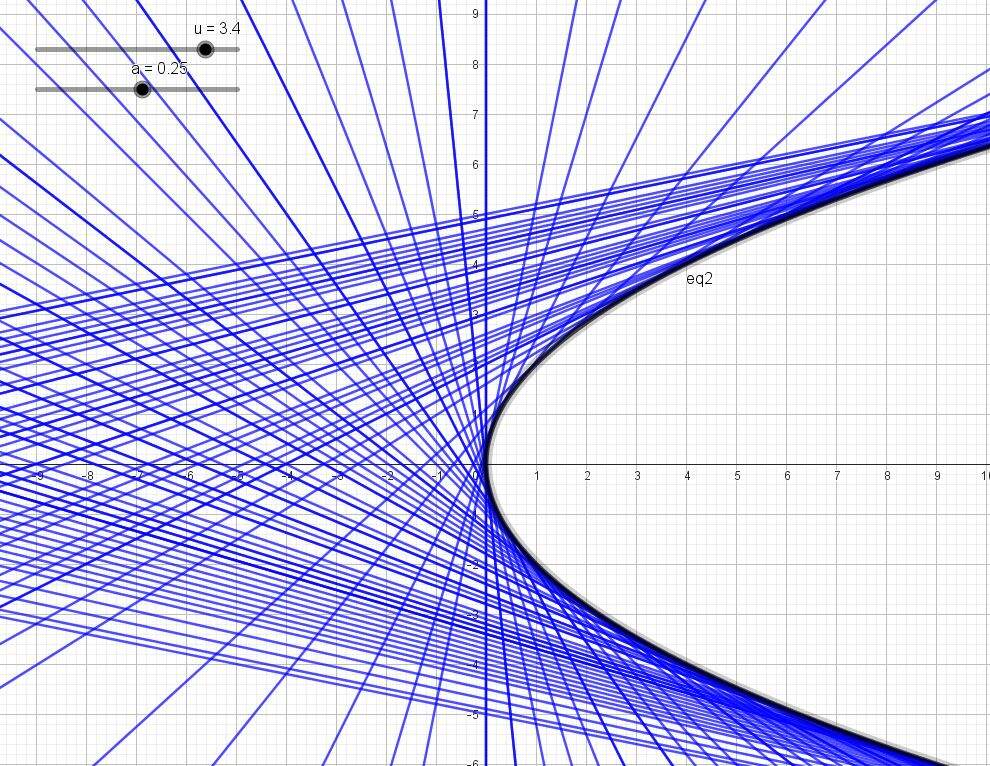,width=4.5cm}}
       }
\caption{Exploration of a family of lines}
\label{fig lines-CAS}
\end{figure}

Figure~\ref{fig lines-CAS} shows a GeoGebra session using its embedded CAS to determine characteristic points, according to Definition~\ref{def envelope synthetic}. These points have coordinates $M_u=(u^2+u\epsilon,-2u-\epsilon)$. We have:
$\underset{\epsilon \rightarrow 0}{\lim} M_{u+\epsilon} = (u^2,-2u)$. This is a parametric presentation of the parabola whose equation is $y=-4x^2$. Figure~\ref{fig lines-CAS}(b) models two neighboring lines $L_u$ and $L_{u+\epsilon}$. Their intersection is given by:
\begin{equation}
\label{eq 2 neighboring lines}
\begin{cases}
x+uy+u^2=0\\
x+(u+\epsilon)y+(u+\epsilon)^2=0
\end{cases}
\end{equation}
The obtained parabola is plotted in Figure~\ref{fig lines-CAS}(c).

Applying this method in the general case shows that
Definition~\ref{def envelope synthetic} implies Definition~\ref{def
  envelope analytic}; see~\cite{bruceandgiblin,dpzrevival}.  GeoGebra
has a command for computing envelopes, but its syntax does not fit the
above problem. The command has syntax Envelope($<$Path$>,<$Point$>$);
this means that the command has been constructed to determine the
envelope of a family of curves (the Path) depending geometrically on a
tracer (the Point). It does not fit a purely algebraic setting, as in
the example above.

\subsubsection{A nephroid as an envelope.}

\begin{figure}[htb]
\centering
\mbox{
\subfigure[Exploration with the mouse]{\epsfig{file=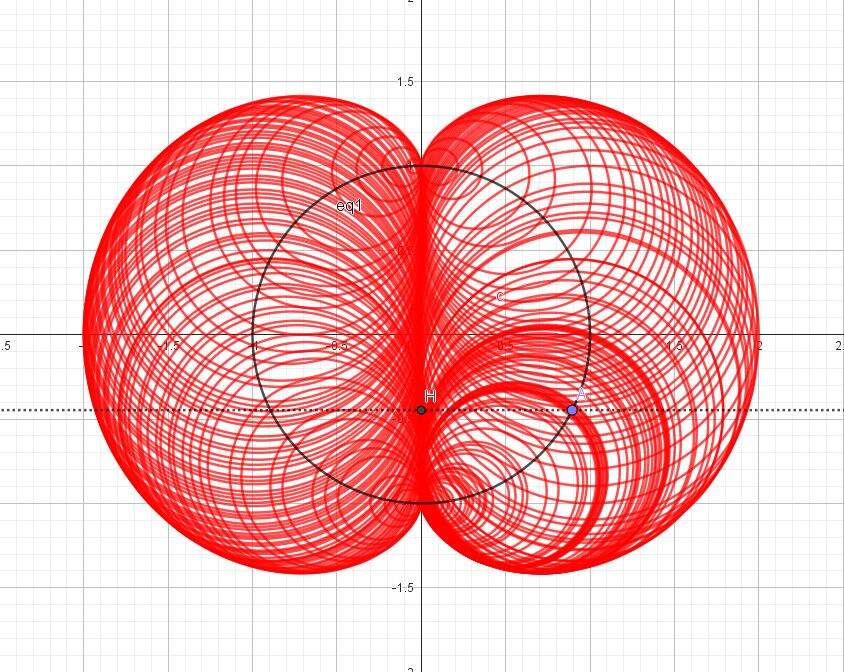,width=6cm}}
\qquad
\subfigure[The Envelope output]{\epsfig{file=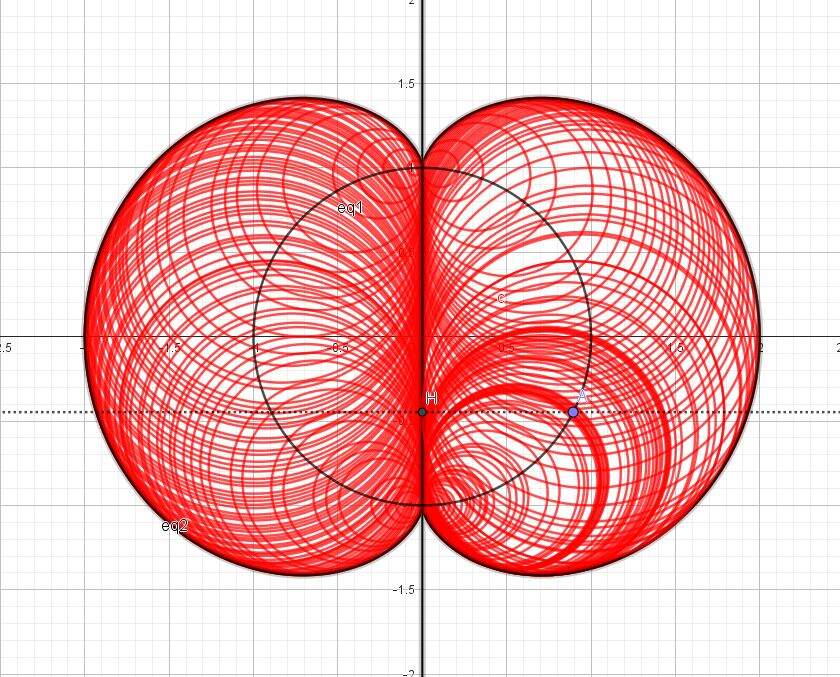,width=6cm}}
       }
\caption{Exploration of a nephroid as envelope of circles}
\label{exploration nephroid}
\end{figure}

Let $\mathcal{U}$ be the unit circle centered at the origin.  Consider the point $A$ on $\mathcal{U}$ and a circle $\mathcal{C}_A$ centered at $A$ and tangent to the $y-$axis. If it exists, denote by $\mathcal{E}$ the envelope of the family of circles. Figure~\ref{exploration nephroid}(a) shows a screen shot of a first experimentation with GeoGebra\footnote{\url{https://www.geogebra.org/m/fxjgmkqu}}, using Trace On for the circles.

Figure~\ref{exploration nephroid}(b) shows a screenshot of the applet after usage of the Envelope command.
\begin{itemize}
\item Draw a line through $A$ perpendicular to the $y-$axis;
\item Determine the point of intersection of this line with the $y-$axis, denoted by $H$ (this may be different with every new experimentation);
\item Plot a circle whose center is $A$ and which passes through $H$.
\end{itemize}
The Envelope command is effective. The output has two components:
\begin{itemize}
\item A curve plotted in the geometric window (the dotted curve of the figure);
\item An implicit equation in the algebraic window.
\end{itemize}
Note that the equation is of degree 7. Actually the polynomial is reducible and can be written as follows:
\begin{equation}
\label{eq nephroid as an envelope}
x(4x^6+12x^4y^2-12x^4+12x^2y^4-24x^2y^2-15x^2+4y^6-12y^4+12y^2-4)=0
\end{equation}
This means that the result is the union of the $y-$axis (given by the vanishing of the first term) and a sextic. A quick websearch yields that this sextic is a nephroid. Note that, if  the Mathcurve site (\url{https://mathcurve.com/courbes2d.gb/nephroid/nephroid.shtml}) is read, the equation given there has to be modified, because of the respective roles of the coordinate axes.
It should be observed that the $y-$axis is ``too big''. It corresponds to the factor $x$ in the polynomial of degree 7  (Equation (\ref{eq nephroid as an envelope}) obtained in the algebraic window, but the geometric data indicates that only a segment of the axis is relevant. The superfluous parts are a consequence of the algebraic computations which determine a closure in the Zariski topology.

The question is now: can we improve the automated work in order to obtain the ``true'' answer? The answer is yes, as we show now.

Consider the following parametric presentation for the unit circle $(x,y)=(\cos u, \sin u), \; u \in [0,2\pi ]$. Then a generic equation for the circles is
\begin{equation}
\label{eq circles - trig}
(x-\cos(u))^2+(y-\sin u)^2-\cos^2u=0.
\end{equation}
Denote by $F(x,y,u)$ the left hand-side in Equation (\ref{eq circles - trig}). In this case, the system of equations of Definition~\ref{def envelope analytic}
 reads as follows:
 \begin{equation}
 \label{system eq nephroid}
\begin{cases}
(x-\cos(u))^2+(y-\sin u)^2-\cos^2u=0 \\
x \sin u - y \cos u +\sin 2u =0
\end{cases}
\end{equation}
Solving the system with Maple, we obtain the following output:

\begin{verbatim}
{x = 0, y = sin(u)}, {x = 2*cos(u)^3, y = -2*sin(u)^3 + 3*sin(u)}
\end{verbatim}

The first component describes a segment on the $y-$axis, which is really accurate. The second component is a parametric presentation of the nephroid. Figure~\ref{animate circles and nephroid} shows two screenshots of an animation performed by Maple:
\begin{figure}[htb]
\centering
\mbox{
\subfigure[]{\epsfig{file=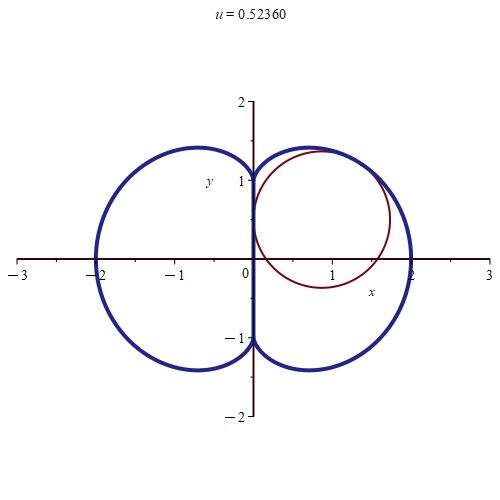, width=5cm}}
\quad
\subfigure[]{\epsfig{file=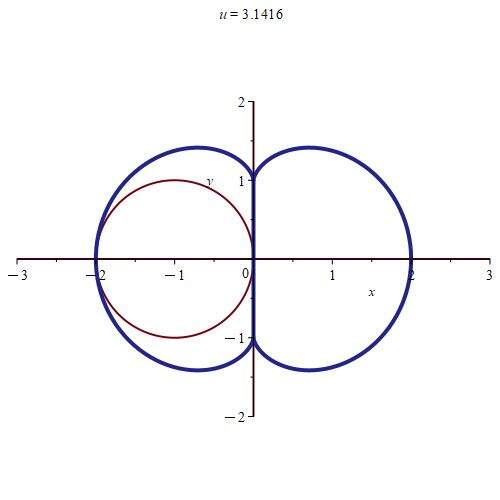, width=5cm}}
     }
\caption{Plots of circles and together with the envelope with Maple}
\label{animate circles and nephroid}
\end{figure}
For the reader's sake, we include here the Maple code. Note that the axes have been defined explicitly in order to have a more readable plot.
\small
\begin{verbatim}
F := (x - cos(u))^2 + (y - sin(u))^2 - cos(u)^2;
derF := diff(F, u);
expand(%);
solve({F = 0, derF = 0}, {x, y});
neph := plot({[0, sin(u), u = 0 .. 2*Pi],
 [2*cos(u)^3, -2*sin(u)^3 + 3*sin(u), u = 0 .. 2*Pi]},
 scaling = constrained, thickness = 4, color = navy);
axes := implicitplot({x = 0, y = 0}, x = -3 .. 3, y = -2 .. 2);
circles := animate(implicitplot, [F = 0], u = 0 .. 2*Pi, thickness = 2);
display(axes, neph, circles);
\end{verbatim}
\normalsize

From this point, the derivation of an implicit equation requires
applying the methods described
in~\cite{talbot,DPKdialog,isopticsfermat,dp-mozgawainner},
transforming trigonometric functions into rational expressions, then
transforming the data into polynomials which generate ideals in a
polynomial ring, and finally utilising elimination algorithms.

\subsection{Isoptics of plane curves.}
\label{subsection isoptics}
\begin{figure}[htb]
\centering
\mbox{
\subfigure[Isoptics of an ellipse]{\epsfig{file=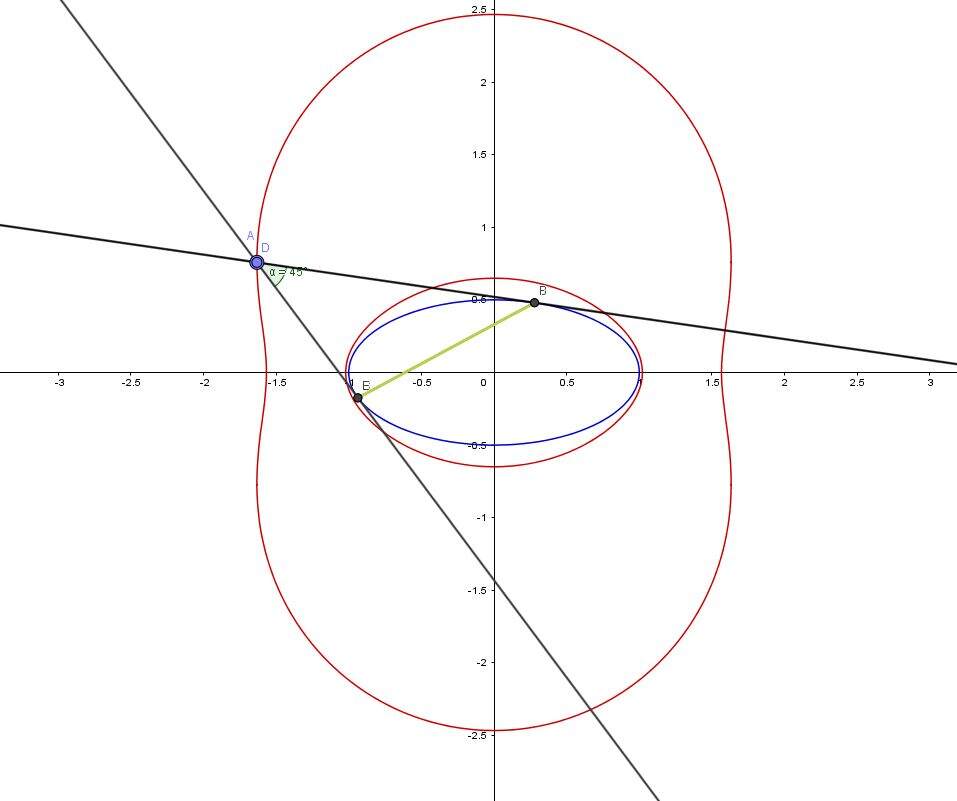, width=6cm}}
\qquad
\subfigure[Isoptic of a Fermat curve]{\epsfig{file=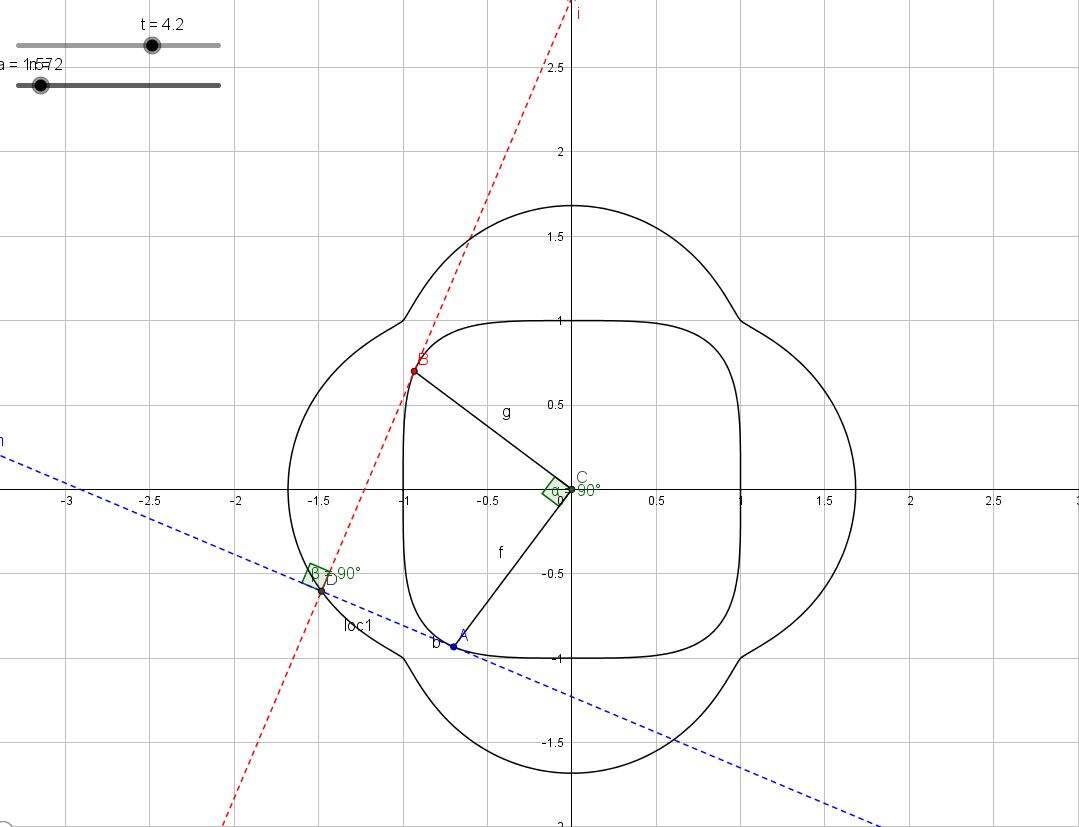,width=6cm}}
     }
\caption{Isoptic curves}
\label{fig isoptics}
\end{figure}

Let $\mathcal{C}$ be a plane curve and $\theta$ a given angle. The
$\theta-$isoptic curve of $\mathcal{C}$ is the geometric locus of
points $M$ in the palen through which passes a pair of tangents making
an angle equal to $\theta$. For conics and $\theta=90^o$, the isoptic
curves are called orthoptics and have been known for a long
time. These are the directrix of a parabola, the director circle of an
ellipse, and the director circle of a hyperbola (when it exists,
depending on the angle between the asymptotes). On Figure~\ref{fig
  isoptics}(a) the isoptics an ellipse are displayed for two
complementary angles. Figure~\ref{fig isoptics} shows the orthoptic
curve (i.e. the isoptic for right angles) of a Fermat curve of degree
6.

For general $\theta$, isoptic curves of conics and of Fermat curves
have been studied, for example,
in~\cite{cieslak-miernowski-mozgawa1991,miernowski-mozgawa1997,DMZ-bisopticsofellipses,DMZ-bisopticsofhyperbolas,isopticsfermat}. The
study relies strongly on the usage of software, but not on
isoptics-dedicated commands. The main tools are solvers for non-linear
equations, and Gr\"obner packages.

In these works, plots of isoptics are provided, but as a consequence of algebraic computations and on a one by one basis. A first step towards automated exploration has been made, by developing tools for automated coloring~\cite{DPKdynamiccoloring,DPK-dynamicalcoloring}. These tools provide visualization together of several isoptics. The variations of the parameter are translated into coloring differences. Confirmation is obtained a posteriori, using GeoGebra's dragging of a point on the given curve\footnote{The interested reader can try the applets \url{https://www.geogebra.org/m/a2zpetsc} for isoptics of an ellipse, \url{https://www.geogebra.org/m/kvbnpzt3} for inner isoptics of an ellipse, and \url{https://www.geogebra.org/m/yjgsbbpk} for isoptics of a Fermat curve. Other applets are also available on the site.}  .

In every case, the availability of the automated command Locus enabled to plot immediately the curve. Moreover, as in other situations,  a slider (or more than one) offer the possibility to explore not only one curve at a time, but a family depending continuously on one (or more) parameter. Nevertheless, this command did not provide an implicit equation of the curve. The derivation of an implicit equation requires algebraic machinery. Here a command to solve equations is needed; as mentioned previously, such  a command is a standard in any CAS. Of course, it has to contain pattern recognition, as different problems lead to different algebraic settings (polynomial equations, trigonometric equations, etc...).

The reward of such activities is multiple:
\begin{itemize}
\item Work is interactive, i.e. it is performed as a dialog between man and machine. This is a nice opportunity to experience the mutual influence between the user and the software, and the teacher can observe a different processes of instrumental genesis among the students; see~\cite{trouche,artigue,artiguetrouche2021}.
\item Generally, the catalogues of curves present them as a discrete set of objects. The automated exploration reveals often that different curves belong actually to one larger family, the various curves corresponding to different values of the parameters. The author described such a situation in~\cite{talbot}, where the machine-and-machine communication is crucial.
\item These activities can stand by themselves, but can also be a
  promo aimed at broadening horizons. When working in class, a couple
  of students reacted with a great ``waooo!'' when the curves appeared.
\end{itemize}

A central issue for the determination of isoptics and envelopes is
solving systems of non linear equations. The result consists in
parametric representation of plane curves. these parametric equations
are rarely rational. Computational work may sometimes been performed
in order to transform the expression into rational expressions and
afterwards to obtain polynomials. If possible, Gr\"obner methods (such
as Elimination) can be applied to implicitize the parametric
presentation. This is not always possible. Even when polynomials can
be obtained, Elimination may be too time-consuming or too
memory-consuming for the CAS to answer; we experienced that when
preparing~\cite{KDR2022}: the computations have been launched on two
different computers with different characteristics, none gave an
answer.

\section{Conclusions and some thoughts for next steps}

We presented activities around automated exploration, discovery and
proof using GeoGebra and Maple. Actually a few drops in a vast
ocean. First, note that we used the current version of GeoGebra, and
the companion package called GeoGebra-Discovery, developed by
Z. Kov\'acs.  It contains several automated commands, sometimes
specific to the package, sometimes an extension of commands existing
in GeoGebra. For example, the Relation command exists in GeoGebra and
is numerical. The GeoGebra-Discovery version extends it with symbolic
algorithms; we showed an example in subsection~\ref{subsection machine
  interactions}.  Other systems are also available for the same kind
of tasks. For example see~\cite{botana-abanades2011} on the usage of
Sage for automated work. An important decision to be made by the
educator is to choose the software which will be used. This choice has
multiple faces. First of all, it depends on which systems are
available at the teacher's institution. This is a consequence of
pedagogical choices, but not less on financial decision of the
administrators. It depends also on the teacher's literacy with regards
to the different available packages. Of course, not every CAS or DGS
is suitable for every level of students. Button-driven software may be
easier to use than a software which requires mastering the syntax of
the commands. Moreover, button-driven commands are often accompanied
by a description of the command when right-clicking on the button, a
very helpful feature.  The existence of an interactive website for
examples and tutorials is a strong help for the student.

The level of the students' background is crucial. Automated methods
are not intended to be used as a blackbox, but are intended to incite
the user to develop new approaches, and to achieve a more profound
insight and understanding of the questions under study. The
\emph{black box - white box} issue has been analyzed for example
in~\cite{drijvers1995,drijvers2000}). Sometimes, the CAS helps to
bypass a lack of theoretical knowledge, but then it is important to
``go back'' and have the students fill the gap and understand what was
hidden in the activities and the software
usage~\cite{bypass}. Critical Thinking and Creativity have to be at
work together. Before presenting the new methods, the teacher has to
be informed of the actual level of the students, and then will be able
to construct a curriculum adapted to the exploration of their Zone of
Proximal Development (ZPD); see~\cite{vygotsky1978}. We recall that
``the ZPD refers to the learner's ability to successfully complete
tasks with the assistance of more capable other people, and for this
reason it is often discussed in relation to assisted or scaffolded
learning. The creation of ZPDs involves assistance with the cognitive
structuring of learning tasks and sensitivity to the learner's current
capabilities'' (Walker,~\cite{walker}). The author of this survey
insists that his students work by pairs. The communication between
them and their subsequent collaboration contribute to a reinforcement
of the benefits that receive from the scaffolding by the teacher. The
student's ``waoow!'' mentioned previously was such a consequence, and
the mathematical discovery made on that occasion was an incitement to
create something new. That has been done.

The versatility of exploratory tools offered by a DGS and a CAS enable
different students in the same class to experience difference ways to
solve the same problem. For STEAM oriented students, this opens
numerous opportunities to create models, animations, etc.  We
illustrated that in subsection~\ref{subsection envelopes}. The
scaffolding offered by the teacher has therefore to be more personally
adapted, and the exploration of the student's ZPD becomes more and
more personal. Some kind of joint ZPD (for 2 students) is also
created. Smartphones, and before that walkmen (who remember them? and
other electronic devices, created some disconnection between
humans. Here, the fact that students are more free to experience and
explore with their personal device, and then share their discoveries
with their classmates, makes the class experience richer, more
meaningful and more interesting. Outdoor activities contribute also to
attract students to learn more mathematics, as these are based on
their everyday environment and to their cultural background. The
author teaches practical courses on ``technology in Mathematics
education'' aimed either at pre-service or in-service teachers. The
technological discourse has always to be adapted, according to the
mathematical level, the cultural background and other characteristics
of the students. Artigue~\cite{artigue} emphasizes also that the new
technological knowledge is an integral part of the new mathematical
knowledge acquired by the students.  New topics or renewed topics can
be proposed and explored, at an earlier stage of the student's cursus
than in the past, whence a need to analyse new didactic situations and
didactic transposition~\cite{brousseau1997,brousseau-warfield}. All
this is part of the new paradigms evoked by
P. Quaresma~\cite{quaresma2020} (v.s. Section~\ref{intro}). It is not
surprising that the instrumental genesis
(see~\cite{trouche,artigue,artiguetrouche2021}, etc.) of every
single student and of the class as a whole is totally different every
year. Of course, this is relevant also for pre-service engineers.

Finally, we should once again emphasize that automated methods for
exploration, discovery and proof, are a new approach aimed at
developing new skills and to emphasize a new kind of
understanding. New developments are always made. Quoting once again
P. Quaresma (op.cit.): ``Geometric reasoning with such computer
applications is one of the most attractive challenges for future
accumulation and dissemination of knowledge.'' For example, recent
works are aimed at the study of inequalities~\cite{inequalities2021}
and their plots in the plane. Something was already available with
Derive, but modern developments are richer. This kind of developments,
together with the personal teaching-learning processes deserve study
of a new kind of instrumental genesis. From a another point of view,
with new achievements for an efficient and fruitful automated dialog
between different kinds of software, new pairs CAS-DGS (either
distinct or embedded one on the other) provide a new artifact which
has to be transformed into an instrument~\cite{artigue,trouche}. A new
loop in a cognitive-educative spiral.


\bibliographystyle{eptcs}

\bibliography{computerAssistedProofsAndAutomatedMethodsInMathEduc17012023.bib}
\end{document}